\documentclass{article}
\usepackage{amsthm,amssymb,amsmath}

\newtheorem{thm}{Theorem}[section]
\newtheorem{lem}[thm]{Lemma}
\newtheorem{prop}[thm]{Proposition}

\newtheorem{rem}[thm]{Remark}

\title{Global existence and decay estimates for the nonlinear wave equations with space-time dependent dissipative term}
\author{TOMONARI WATANABE\footnote{Hiroshima University, e-meil:watanabe-to@hiroshima-u.ac.jp}}

\begin{document}

\maketitle
 
\begin{abstract}
We study the global existence and decay estimates for nonlinear wave equations with the space-time dependent dissipative term in an exterior domain.
The linear dissipative effect may vanish in a compact space region. Moreover the nonlinear terms need not divergence form.
For getting the higher order energy estimates, we introduce an argument using the rescaling. The method is useful to control derivatives of the dissipative coefficient.


\end{abstract}

\numberwithin{equation}{section}
\allowdisplaybreaks

\section{Introduction}
Let $d\geq2$ and $\Omega = \mathbb{R}^d/ \mathcal{O}$, where $\mathcal{O}$ is a star-shaped domain with a smooth and compact boundary $\partial \Omega$.
Moreover we assume that $\mathcal{O}$ contains the origin.
In this paper, we consider the initial-boundary value problem for nonlinear wave equations with the space-time dependent dissipative term:
\begin{equation}\nonumber
\rm{(DW)}  \hspace{2mm} \left\{
\begin{array}{ll}
\displaystyle (\partial_t^2 - \triangle + B(t,x)\partial_t) u (t,x) = F(\partial u, \partial^2 u)  & (t,x) \in [0,\infty) \times \Omega, \\
u(0,x) = u_0(x) ,\quad  \partial_t u(0,x) = u_1 (x)  & x \in \Omega,\\
u(t,x) = 0  & (t,x) \in [0,\infty) \times \partial \Omega,
\end{array}
\right.
\end{equation}
where $u=(u^1, \cdots, u^d)$, $\nabla = (\partial_{x_1}, \cdots, \partial_{x_d})$ and $\partial = (\partial_t, \nabla)$. 
The initial data $(u_0,u_1)$ belongs to $H^{L}(\Omega) \times H^{L-1}(\Omega)$ and satisfies the compatibility condition of order $L-1$. $H^L(\Omega)$ is the Sobolev space in $\Omega$.
We make the following assumptions for the space-time dependent damping coefficient matrix $B(t,x)=(B_{pq}(t,x))_{p,q = 1,\cdots,d}$:
\begin{description}
\item[(B0)] $B_{pq}$ belong to $\mathcal{B}^\infty ([0,\infty) \times \Omega)$, where $\mathcal{B}^\infty$ is the function space of smooth functions with bounded derivatives.
\item[(B1)] $B(t,x)$ is nonnegative definite in $[0,\infty) \times \Omega$.
\item[(B2)] $\partial_t B (t,x)$ is nonpositive definite in $[0,\infty) \times \Omega$.
\item[(B3)] There exist $b_0 >0$ and $R>0$ such that 
\begin{equation}\nonumber
\sum_{p,q=1}^dB_{pq}(t,x)\eta_p \eta_q \geq b_0 |\eta|^2 \quad (t \in [0,\infty), |x| \geq R, \eta \in \mathbb{R}^d ).
\end{equation}
\end{description}
As is in {\bf (B3)}, dissipative term works on $|x| \geq R$. This means that dissipative effect may vanish in a compact space region. 

We treat quadratic nonlinear terms. 
In what follows, $\partial_0$ means $\partial_t$ and $\partial_j (j=1,2, \cdots, d)$ means $\partial_{x_j}$.
Assume that $F$ is of the form 
\begin{equation}\nonumber
\displaystyle F(\partial u,\partial^2 u)= \left(\tilde{F}_i( \partial u) + \sum_{j=1}^d \sum_{0 \leq a,b \leq d} c_{ij}^{ab}(\partial u) \partial_a \partial_b u^j \right)_{i=1,\cdots,d},
\end{equation}
which satisfy 
\begin{equation}\label{eq:condi0}
c_{ij}^{ab} = c_{ji}^{ba},
\end{equation}
\begin{equation}\label{eq:condi1}
|D_\xi^\alpha \tilde{F}_i(\xi)|\leq C_{\alpha,p_1} |\xi|^{\max\{0, p_1-|\alpha|\}} \quad  ( \xi \in \mathbb{R}^d \times \mathbb{R}^{d+1}, \ \  |\alpha| \leq L-1)
\end{equation}
and
\begin{equation}\label{eq:condi2}
|D_\xi^\alpha c_{ij}^{ab} (\xi) | \leq C_{\alpha,p_2} |\xi|^{\max\{0, p_2-1-|\alpha|\}} \quad (\xi \in \mathbb{R}^d \times \mathbb{R}^{d+1}, \ \  |\alpha| \leq L-1)
\end{equation}
for some $p_l\geq 2\ (l=1,2)$.
The main objective of this paper is to prove the global existence and decay estimate to (DW).

In the case of the coefficient function $B$ vanishes, (DW) become the nonlinear wave equations.
Then it is well known that no matter how small the initial data, there do not exist globally defined smooth solutions in general (e.g.\cite{John1}, \cite{Kubo}, \cite{Keel}).
If $F$ has the "Null condition" then (DW) has a global smooth solution for sufficiently smooth and small the initial data (e.g.\cite{Klainerman}, \cite{Sideris}).

In the case of the coefficient function $B \equiv Const >0$, there are many results (\cite{Kawashima}, \cite{Ma1} etc.).
When linear or semilinear version, it is known that the asymptotic profile of the solution to (DW) is given by the corresponding solution of heat equation (e.g. \cite{Narazaki}, \cite{Nishihara} etc.).
Such a property is called the diffusion phenomenon.
Recently, there are many research concerning the diffusion phenomenon for the nonuniform dissipation.
When the dissipation depends on space variable $B=B(x)$, Todorova and Yordanov consider like $B= (1+|x|)^{-\gamma}$ to linear and semilinear version in \cite{Todorova}, \cite{Todorova2}.
They show that if $0 \leq \gamma < 1$ then the solution to (DW) have some decay estimates, indeed Wakasugi \cite{Wakasugi} confirms the diffusion phenomenon recently.
When the dissipation depend on time variable $B=B(t)$, Wirth \cite{Wirth} proves that if $tB(t) \to +\infty$ and $B \in L^1$ then the solution to (DW) satisfy the decay estimate like corresponding solution to heat equation.
On the other hand, Mochizuki \cite{Mochizuki} consider the scattering for the free wave equation when $B$ depends time-space variable.
He prove that if there exist $\xi,\eta \in L^1$ and small $\varepsilon$ such that
\begin{equation}\nonumber
|B(t,x)| \leq \varepsilon \xi(|x|) + \eta(t),\ \ \xi, \eta \geq 0,\ \xi' \leq 0, \ \  \xi'^2 \leq 2 \xi \xi''
\end{equation}
then the solution to (DW) close to the free wave equation.
We remark that $B(t,x) = (1+t)^{-\alpha} (1+|x|)^{-\beta},\ \alpha + \beta >1$ is a sufficient condition of above the conditions.

Now we consider the nonuniform dissipative term which doesn't decay near infinity but vanishes in a compact region.
Nakao \cite{Nakao} gets the energy decay estimates like $E(u(t)) = O((1+t)^{-1})$ when $B$ depends on space variable only, where $E(u(t))$ is the standard energy of wave equations.
Furthermore, Ikehata \cite{Ikehata2} get the decay estimates as $\|u(t)\|_{L^2(\Omega)}^2 = O((1+t)^{-1})$ and $E(u(t)) = O((1+t)^{-2})$ with the additional condition for the initial data. 
Those results mostly deal with the linear and the semilinear problem. 
For quasilinear version with divergence form, Nakao \cite{Nakao} consider the equation 
\begin{equation}
\partial_t^2u - {\rm div}\{\sigma(|\nabla^2u|) \nabla u\} + a(x) \partial_tu=0,
\end{equation}
where $\sigma$ is a smooth function like $\sigma(x) = (1+|x|)^{-\frac{1}{2}}$. 
Then the nonlinearity order $p$ satisfy $p \geq 3$. Besides, the author of the present paper deals with $p=2$ and prove the decay estimates in \cite{W}.
Moreover there is no result when $B$ depend space-time variables.
In this paper, we assume the space-time dependent dissipation $B$ effective near the space infinity, even if the nonlinear terms $F$ have no "null condition" and "divergence form".

First, we get the global existence as follows:
\begin{thm}\label{thm:main}
Let $L \geq [\frac{d}{2}] + 3$.
Then there exists a small constant $\hat{\delta}>0$ such that if the initial data $(u_0,u_1) \in H^{L}(\Omega) \times H^{L-1}(\Omega)$ satisfy the compatibility condition of order $L-1$ and
\begin{equation}\label{eq:1.4}
\|(u_0,u_1) \|_{H^{L}(\Omega) \times H^{L-1}(\Omega)} \leq \delta, 
\end{equation}
then there exists a unique global solution to {\rm (DW)} in $\bigcap_{j=0}^{L-1} C^j([0,\infty); H^{L-j}(\Omega) \cap H^1_0 (\Omega) ) \bigcap C^L([0,\infty); L^2 (\Omega))$.
\end{thm}
In the proof of Theorem \ref{thm:main}, we use higher order energies (see for instance \cite{Kubo},\cite{Sideris}) and the rescaling (see section 2). 
Note that if $B=Const>0$, we can prove theorem \ref{thm:main} under the assumption $\|(\nabla u_0,u_1)\|_{H^{L-1(\Omega)} \times H^{L-1}(\Omega)} \leq \hat{\delta}$ instead of \eqref{eq:1.4}, 
i.e. the smallness of $\|u_0\|_{L^2(\Omega)}$ is also needed for the case of nonuniform dissipative term.
 
Next, we also prove the decay estimate as follows:
\begin{thm}\label{thm:main2}
In addition to the assumptions in Theorem \ref{thm:main}, we assume {\rm(H1)} and {\rm(H2)}.
\begin{description}
 \item[(H1)] $\| d_0(\cdot) \{ B(0)u_0 + u_1 \} \|_{L^2(\Omega)} <\infty $.
 \item[(H2)] $\displaystyle  \int_0^\infty \| d_0(\cdot) \partial_t B (s) \|_{L^\infty(\Omega)} ds < \infty$,
\end{description}
where $d_0:\mathbb{R}^d \to \mathbb{R}$ is defined by

\begin{equation}\label{eq:d_0}
d_0(x) =\left\{
\begin{array}{ll}
|x|   &\quad (d \geq 3), \\
|x| \log(A |x|)  &\quad (d = 2)
\end{array}
\right.
\end{equation}
with a constant A satisfying $\inf_{x \in \Omega } A |x| \geq 2$.
Furthermore if $d=2$, we also assume {\rm(H3)}.
\begin{description}
 \item[(H3)] There exists $M$ such that ${\rm supp} u_0 \cup {\rm supp} u_1 \subset \{x \in \Omega: |x| \leq M \}$.
\end{description}
Then the global solution $u$ to {\rm (DW)} satisfy following estimates:
\begin{equation}\label{eq:1.6}
\sum_{\mu=0}^{L-1} \|(\partial_t^\mu u(t),\partial_t^{\mu+1} u(t))\|^2_{H^{L-\mu}(\Omega) \times H^{L-\mu-1}(\Omega)} \leq E_0(1+t)^{-1}
\end{equation}
and
\begin{equation}\label{eq:1.7}
\|(\nabla u(t),\partial_t u(t))\|^2_{L^2 (\Omega) \times L^2 (\Omega)} \leq E_0(1+t)^{-2},
\end{equation}
where $E_0$ is a constant depend on $(u_0,u_1)$ and $M$.
\end{thm}
Theorem \ref{thm:main2} is the decay estimate which correspond to the result of Ikehata \cite{Ikehata2}. 

The paper is organized as follows.
In section 2 we prepare some known lemmas and the rescaling function. 
In section 3 we prove the high energy estimates to (DW) and the theorem \ref{thm:main}. In the proof, the rescaling argument plays an important role.
In section 4 we prove the Theorem \ref{thm:main2}.

\section{Preliminaries}

We consider the rescaling to (DW).
Let $u$ be the solution to (DW). We define $\displaystyle v(t,x) = \lambda^{-1} u\left(\lambda t, \lambda x\right) (\lambda>0)$, then $v$ satisfies
\begin{eqnarray*}
&&\partial_t^2 v (t,x)  - \triangle v (t,x) = \lambda\left\{ \partial_t^2 u\left(\lambda t, \lambda x\right) - \triangle u\left(\lambda t, \lambda x\right)  \right\}  \\
 &=& -\lambda B \left(\lambda t,\lambda x\right) \partial_t u \left(\lambda t, \lambda x\right) + \lambda F( \partial u(\lambda t, \lambda x), \partial \nabla u(\lambda t, \lambda x))  \\
 &=& -\lambda B \left(\lambda t,\lambda x\right) \partial_t v(t,x) + F_\lambda(\partial v(t,x), \partial \nabla v(t,x)),
\end{eqnarray*}
where $\displaystyle (F_{\lambda})_i (\partial v , \partial \nabla v ) = (\tilde{F}_{\lambda})_i(\partial v) + \sum_{j=1}^d \sum_{0 \leq a,b \leq d} c_{ij}^{ab}(\partial v) \partial_a \partial_b v_j $ and $\tilde{F}_{\lambda}(\partial v) = \lambda \tilde{F}( \partial v)$.

So $v$ is the solution to the following initial-boundary value problem (DW$)_\lambda$:
\begin{equation}\nonumber
\rm{(DW)_\lambda}  \hspace{2mm} \left\{
\begin{array}{ll}
\displaystyle (\partial_t^2 - \triangle + B_\lambda(t,x)\partial_t ) v = F_\lambda (\partial v,\partial \nabla v)  &  [0,\infty) \times \Omega_\lambda , \\
v(0,x) = v_0(x) ,\quad  \partial_t v(0,x) = v_1 (x)  & x \in \Omega_\lambda, \\ 
v(t,x) = 0 & (t,x) \in [0,\infty) \times \partial \Omega_\lambda,
\end{array}
\right.
\end{equation}
where $\Omega_\lambda = \{x : \lambda x\in\Omega \}$, $ B_\lambda(t,x) = \lambda B\left(\lambda t, \lambda x \right) , v_0 (x)= \lambda^{-1} u_0(\lambda x) , v_1(x) = u_1(\lambda x)$.
Then $B_\lambda$ satisfy following \textbf{(B1$)_\lambda$} - \textbf{(B3$)_\lambda$} instead of \textbf{(B1)} - \textbf{(B3)}. 
\begin{description}
\item[(B1$)_\lambda$] $B_\lambda(t,x)$ is nonnegative in $[0,\infty) \times \Omega_\lambda$,
\item[(B2$)_\lambda$] $\partial_t B_\lambda(t,x)$ is nonpositive in $[0,\infty) \times \Omega_\lambda$,
\item[(B3$)_\lambda$] There exist $b_0 >0$ and $R>0$ such that 
\begin{equation}\nonumber
\sum_{p,q = 1, \cdots, d}(B_\lambda(t,x))_{pq} \eta_p \eta_q  \geq \lambda b_0 |\eta|^2 \qquad (t \in [0,\infty) , |x| \geq \frac{R}{\lambda}, \eta \in \mathbb{R}^d ).
\end{equation}
\end{description}
Furthermore $B_\lambda$ satisfies
\begin{description}
\item[$\textbf{(B4)}_\lambda$]$ \|\partial^\alpha  B_\lambda\|_{L^\infty([0,\infty) \times \Omega_\lambda)} \leq \lambda^{|\alpha|+1} \|\partial^\alpha B\|_{L^\infty([0,\infty) \times \Omega)} $.
\end{description}
We consider (DW$)_\lambda$ instead of (DW). Throughout this paper, $\| \cdot \|_p = \| \cdot \|_{L^p(\Omega_\lambda)}$, $\|\cdot \|_{H^l} = \|\cdot\|_{H^l(\Omega_\lambda)}$ and $\langle\cdot, \cdot \rangle$ stands for $L^2(\Omega_\lambda)$-inner product.

First, we introduce some known results. First we prepare following lemmas for estimating nonlinear terms.
\begin{lem}[Sobolev's lemma]\label{prop:sobolev}
There exists a constant $C_\lambda> 0$ such that
\begin{equation}\nonumber
\|f\|_\infty \leq C_\lambda \|f\|_{H^{\left[\frac{d}{2}\right] +1}} \quad ( f \in H^{\left[\frac{d}{2}\right]+1}(\Omega_\lambda) ).
\end{equation}
\end{lem}

\begin{lem}[Elliptic estimate]\label{lem:elliptic}
There exists $C_\lambda>0$ such that for any $\varphi \in H^m (\Omega_\lambda) \cap H^1_0 (\Omega_\lambda)$ with an integer $m \geq2$, we have
\begin{equation}\nonumber
\sum_{|\alpha| = m} \|\nabla^\alpha \varphi\|_2 \leq C_\lambda (\|\triangle \varphi \|_{H^{m-2}} + \|\nabla \varphi\|_2).
\end{equation}
\end{lem}

Next, we prepare the Poincare type inequality associated with $B_{\lambda}$.
\begin{lem}[Poincare type inequality]\label{lem:Poincare}
There exists a constant $C_1\geq 1/4$ such that 
\begin{equation}\label{eq:2.1}
\|f\|^2_2 \leq \lambda^{-1}C_1 \langle f,B_\lambda f \rangle + \lambda^{-2} C_1 \|\nabla f\|^2_2 \quad (f \in H^1_0(\Omega_\lambda), \lambda>0)
\end{equation}

\proof
We define $U_r = \{ x \in \Omega_\lambda:|x| \leq r \}$. Using Poincare inequality,  
we obtain the following estimate:
\begin{equation}\nonumber
\int_{U_r} |f(x)|^2 dx \leq r^2 \int_{U_r} |\nabla f(x)|^2 dx \quad (f\in H^{1}_0 (U_r)).
\end{equation}
Let $f\in H^1_0(\Omega_\lambda)$ and $\rho \in C_0^\infty(\mathbb{R}^d) $ be a function satisfying $0 \leq \rho\leq 1, \rho (x) = 1 (|x| \leq 1) , \rho(x) = 0 (|x| \geq \frac{3}{2})$ . 
We define $\rho_\lambda (x) = \rho (\frac{\lambda x}{R})$ for $\lambda >0$, then because of $\rho_\lambda f \in H^1_0(U_\frac{2R}{\lambda})$ we have 
\begin{eqnarray*}
 \|f\|^2_2 &=& \int_{\Omega_\lambda} |\rho_\lambda(x) f(x)|^2 dx + \int_{\Omega_\lambda} (1-|\rho_\lambda(x)|^2) |f(x)|^2 dx \\
 &=& \int_{U_{\frac{2R}{\lambda}}} |\rho_\lambda(x) f(x)|^2 dx + \int_{|x| \geq \frac{R}{\lambda}} (1-|\rho_\lambda(x)|^2) |f(x)|^2 dx \\
 &\leq& \frac{4R^2}{\lambda^2}\int_{U_{\frac{2R}{\lambda}}} |\nabla(\rho_\lambda(x) f(x))|^2 dx +  \int_{|x| \geq \frac{R}{\lambda}} |f(x)|^2 dx \\
 &\leq& \frac{4R^2}{\lambda^2}\int_{\frac{R}{\lambda} \leq |x| \leq \frac{2R}{\lambda}} |\nabla \rho_\lambda(x)|^2 | f(x)|^2 dx + \frac{4R^2}{\lambda^2}\int_{U_{\frac{2R}{\lambda}}} |\rho_\lambda(x)|^2 |\nabla f(x)|^2 dx \\
 &&  + \int_{|x| \geq \frac{R}{\lambda}} |f(x)|^2 dx \\
 &\leq& \left( 4 R^2 \|\nabla \rho  \|_\infty^2 + 1 \right)  \int_{|x| \geq \frac{R}{\lambda}} |f(x)|^2 dx + \frac{4R^2}{\lambda^2}\int_{\Omega_\lambda}|\nabla f(x)|^2 dx  \\
 &\leq& \frac{ 4 R^2 \|\nabla \rho  \|_\infty^2 + 1 }{b_0 \lambda} \langle f,B_\lambda f \rangle + \frac{4R^2}{\lambda^2}\|\nabla f\|^2_2 .
\end{eqnarray*}
Hence we get (\ref{eq:2.1}).
\qed
\end{lem}

Finally, we introduce Hardy inequality and Gagliardo-Nirenberg inequality. We need these in the proof of Theorem \ref{thm:main2} in section 4.
\begin{lem}[Hardy inequality; see for instance \cite{Dan}]\label{lem:hardy}
Let $d \geq 2$. There exists a constant  $C_\lambda > 0$ such that any $f \in H^1_0(\Omega_\lambda)$ satisfies
\begin{equation}\nonumber
\left\|\frac{f}{d_0}\right\|_2 \leq C_\lambda \|\nabla f\|_2,
\end{equation}
where $d_0$ is defined by \eqref{eq:d_0}.
\end{lem}
\begin{lem}[Gagliardo-Nirenberg inequality]\label{lem:GN}
Assume $1 \leq q < d $ and  $\frac{1}{r} = \frac{1}{q}-\frac{1}{d}$.
Then there exists a constant $C_\lambda>0$ such that 
\begin{equation}\nonumber
\|f\|_{r} \leq C_\lambda \|\nabla f\|_{q} \quad (f \in C^\infty_0(\Omega_\lambda)) .
\end{equation}
\end{lem}

\section{Energy estimates}
In this section, we give a proof of the high order energy estimates.

First we introduce some notations.
The energy $E(v(t))$ and the higher order energies $Z_m(v(t)), Z(v(t))$ are defined by
\begin{equation}\label{eq:E}
E(v(t)) = \frac{1}{2} \{\|\partial_t v(t) \|^2_2 +\|\nabla v(t)\|^2_2 \},
\end{equation}
\begin{equation}\label{eq:Zm}
Z_{m}(v(t)) = \sum_{\mu=0}^{L-1-m} \left\{ \|\nabla \partial_t^\mu v(t)\|_{H^m}^2 + \|\partial_t^{\mu+1} v(t)\|^2_{H^m} \right\} \quad (0 \leq m \leq L-1)
\end{equation}
and
\begin{equation}\label{eq:Z}
Z(v(t)) = \sum_{m=0}^{L-1} Z_m(v(t)).
\end{equation}
We sometimes omit $t$ or $v(t)$. Note that $\|v\|_2^2 + Z_{L-1}(v) = \|(v,\partial_t v)\|_{H^L \times H^{L-1}}^2$ holds.
The following result concerning local existence is standard:
\begin{prop}[Local existence; for instance \cite{Kato} or \cite{Nakao}]\label{prop:3.1}
Let $L \geq \left[\frac{d}{2}\right] + 3$ and assume that $\partial \Omega_\lambda$ is smooth.
Furthermore we assume the initial date $(v_0,v_1) \in H^{L}(\Omega_\lambda) \times H^{L-1}(\Omega_\lambda)$ satisfies the compatibility condition of order $L-1$ associated with the {\rm (DW$)_\lambda$}.
Then there exists a unique local solution to {\rm (DW$)_\lambda$} in  
\begin{equation}\nonumber
X^T :=\bigcap_{j=0}^{L-1} C^j([0,T); H^{L-j}(\Omega_\lambda) \cap H^1_0(\Omega_\lambda) ) \bigcap C^L([0,T); L^2(\Omega_\lambda)),
\end{equation}
where $T$ depend on $\|(v_0,v_1) \|_{H^L \times H^{L-1}}$.
\end{prop}
We define the function spaces $X_\delta^T$ and $X_\delta$ as follows:
\begin{equation}\nonumber
X_{\delta}^T = \left\{ v \in X^T :\| v(t) \|_2^2 + Z(v(t)) \leq \delta^2 \quad (0 \leq t \leq T) \right\}
\end{equation}
and
\begin{equation}\nonumber
X_\delta = \left\{v \in X^\infty :\| v(t) \|_2^2 + Z(v(t)) \leq \delta^2 \quad (0 \leq t < \infty) \right\}.
\end{equation}
This section, we prove the following proposition:
\begin{prop}\label{prop:3.2}
There exist $0<\lambda <1$ and $\delta = \delta(\lambda)$ such that the local solution $v \in X_\delta^T$ to {\rm(DW$)_\lambda$} satisfies 
\begin{equation}\label{eq:prop3.2}
\|v(t)\|^2_2 + Z(v(t)) + \int_0^t Z(v(s)) ds \leq C_\lambda \|(v_0,v_1)\|_{H^L \times H^{L-1}}^2.
\end{equation}
\end{prop}

If Proposition \ref{prop:3.2} is holds, then we can prove Theorem \ref{thm:main}.
Indeed, using Proposition \ref{prop:3.1} and Proposition \ref{prop:3.2}, we can prove the unique global existence theorem to (DW$)_\lambda$ by standard continuation argument. 
Furthermore we put $ u(t,x) :=\lambda v(\lambda^{-1}t, \lambda^{-1}x)$, it is easy to see $u$ satisfies the statement of Theorem \ref{thm:main}.

In order to show Proposition \ref{prop:3.2}, we only need the estimates of $Z_0$. Indeed we can prove the next lemma (see for instance \cite{Kubo}.).

\begin{lem}\label{lem:3.3}
For any $\lambda>0$ there exist $\delta = \delta(\lambda)$ and $C_\lambda > 0$ such that a local solution $v \in X_\delta^T$ to {\rm (DW$)_\lambda$} satisfy the following estimates:
\begin{equation}\label{eq:3.5}
Z(v(t)) \leq C_\lambda Z_0 (v(t)) \quad (t \in [0,T))
\end{equation}
and
\begin{equation}\label{eq:3.6}
Z(v(t)) \leq C_\lambda Z_{L-1} (v(t)) \quad (t \in [0,T)).
\end{equation}
\proof
Let $v \in X_{\delta}^T$ is a solution to (DW$)_\lambda$. First, we prove \eqref{eq:3.5}. For $1 \leq m \leq L-1$, it hold that
\begin{align*}
Z_m &= \sum_{\mu=0}^{L-1-m} \left\{ \|\nabla \partial_t^\mu v\|_2^2 + \sum_{2 \leq |a| \leq m+1} \|\partial_t^\mu \nabla^a v\|^2_2 + \|\partial_t^{\mu+1} v\|^2_{H^m} \right\} \\
&\leq  C\left (Z_0 + \sum_{\mu=0}^{L-1-m} \sum_{2 \leq |a| \leq m+1} \|\partial_t^\mu \nabla^a v\|^2_2 +  Z_{m-1} \right). 
\end{align*}
Using Lemma \ref{lem:elliptic} and (DW$)_\lambda$, we get 
\begin{align*}
& \sum_{\mu=0}^{L-1-m} \sum_{2 \leq |a| \leq m+1} \|\partial_t^\mu \nabla^a v\|^2_2 \leq C_\lambda \sum_{\mu=0}^{L-1-m} \left(  \|\partial_t^\mu \triangle v\|^2_{H^{m-1}} + \| \partial_t^\mu \nabla v \|_2^2 \right) \\
&\leq C_\lambda \sum_{\mu=0}^{L-1-m} \left( \|\partial_t^{\mu+2} v\|^2_{H^{m-1}} +  \|\partial_t^\mu( B_\lambda \partial_t v )\|^2_{H^{m-1}} + \| \partial_t^\mu F_\lambda \|_{H^{m-1}}^2 + \| \partial_t^\mu \nabla v \|_2^2 \right) .
\end{align*}
It is easy to see that 
\begin{equation}\nonumber
\sum_{\mu=0}^{L-1-m}  \|\partial_t^{\mu+2} v\|^2_{H^{m-1}} \leq  Z_{m-1},\ \ \sum_{\mu=0}^{L-1-m} \| \partial_t^\mu \nabla v \|_2^2 \leq  Z_{0} 
\end{equation}
and
\begin{align*}
&\sum_{\mu=0}^{L-1-m}  \|\partial_t^\mu( B_\lambda \partial_t v )\|^2_{H^{m-1}} \leq C_\lambda Z_{m-1}.
\end{align*}
Furthermore applying Lemma \ref{lem:A1} and Lemma \ref{lem:A2} to the nonlinear terms, we get
\begin{align*}
Z_m(v(t)) &\leq C_\lambda Z_0(v(t)) + C_\lambda Z_{m-1}(v(t))) + C_\lambda Z(v(t))^2  \\
&\leq C_\lambda Z_0(v(t)) + C_\lambda Z_{m-1}(v(t))) + C_\lambda \delta^2 Z(v(t)) .
\end{align*}
Therefore we obtain inductively that
\begin{equation}\nonumber
Z(v(t)) = \sum_{m=0}^{L-1} Z_m(v(t)) \leq C_\lambda Z_0(v(t)) +C_\lambda \delta^2 Z(v(t))  .
\end{equation}
Choosing $\delta$ small enough depending on $\lambda$, we get \eqref{eq:3.5}.

Next, we prove \eqref{eq:3.6}. Using the same way as for the proof of \eqref{eq:3.5}, for $0 \leq m \leq L-3$, we have
\begin{align*}
Z_m &= \sum_{\mu=0}^{L-1-m} \left\{ \| \nabla \partial_t^\mu v \|_{H^m}^2 + \|\partial_t^{\mu+1} v \|^2_{H^m} \right\} \\
 &=\sum_{\mu=0}^1\| \nabla \partial_t^\mu v\|^2_{H^m} + \sum_{\mu=2}^{L-1-m} \|\nabla \partial_t^\mu v \|^2_{H^m} + \| \partial_t v\|^2_{H^m} + \sum_{\mu=1}^{L-1-m} \| \partial_t^{\mu+1} v \|^2_{H^m}  \\
 &\leq C Z_{L-1}(v) \\
 & + C \sum_{\mu=2}^{L-1-m} \left( \|\nabla \partial_t^{\mu-2} \triangle v \|^2_{H^m} + \| \nabla \partial_t^{\mu-2} (B_\lambda \partial_t v) \|^2_{H^m} + \|\nabla \partial_t^{\mu-2} F_\lambda \|^2_{H^m} \right) \\
 & + C\sum_{\mu=1}^{L-1-m} \left( \|\partial_t^{\mu-1} \triangle v \|^2_{H^m} + \|\partial_t^{\mu-1} (B_\lambda \partial_t v) \|^2_{H^m} + \|\partial_t^{\mu-1} F_\lambda \|^2_{H^m} \right) \\
 &\leq C_\lambda \left( Z_{L-1} + Z_{m+1} + Z_{m+2} + \delta^2 Z\right).
\end{align*}
Similarly we can get $\displaystyle Z_{L-2} \leq C_\lambda Z_{L-1} + C_\lambda \delta^2 Z$, thus it follows that 
\begin{equation}\nonumber
Z(v(t)) = \sum_{m=0}^{L-1}Z_m (v(t)) \leq C_\lambda Z_{L-1}(v(t)) + C_\lambda \delta^2 Z(v(t)).
\end{equation}
Choosing $\delta$ small enough depend on $\lambda$, we get \eqref{eq:3.6}.
This completes the proof of Lemma \ref{lem:3.3}.
\qed
\end{lem}
We consider the estimates of $Z_0(v(t))$.
\begin{lem}\label{lem:3.4}
Let $\mu \leq L-1$ and $\lambda \leq 1$. Then there exists a constant $C>0$ such that a local solution $v$ to {\rm(DW$)_\lambda$} satisfy
\begin{equation}\label{eq:A}
\frac{d}{dt} E(\partial_t^\mu v) + \langle \partial_t^{\mu+1}  v, B_\lambda \partial_t^{\mu+1} v \rangle  \leq C \lambda^2 Z_0(v) + \langle \partial_t^{\mu+1}  v , \partial_t^\mu F_\lambda \rangle,
\end{equation}
\begin{align}\label{eq:B}
&\frac{d}{dt} \{ \langle \partial_t^\mu v, \partial_t^{\mu+1} v \rangle + \frac{1}{2} \langle \partial_t^\mu v, B_\lambda \partial_t^\mu v \rangle \} + \|\nabla \partial_t^\mu v \|^2_2 - \| \partial_t^{\mu+1} v\|^2_2 \\
& \hspace{30mm} \leq  C\lambda^2 Z_0(v) + \langle \partial_t^\mu v , \partial_t^\mu F_\lambda \rangle \nonumber
\end{align}
and for any $K>0$, it holds that
\begin{align}\label{eq:C}
&\frac{d}{dt} \langle \partial_t^{\mu+1} v,[h ; \nabla \partial_t^\mu v ] \rangle \\
 &+ \int_{\Omega_\lambda} \left( \frac{d \phi + |x| \phi'}{2} \right) |\partial_t^{\mu+1}v|^2 dx + \int_{\Omega_\lambda} \left( \frac{(2-d)\phi+|x|\phi'}{2}  \right)|\nabla \partial_t^\mu v|^2  dx \nonumber \\ 
 & \leq \frac{1}{2} \int_{\partial \Omega_\lambda} h \cdot \sigma | \sigma \cdot \nabla \partial_t^\mu v|^2  dS+ \frac{K}{4}\langle \partial_t^{\mu+1} v,B_\lambda \partial_t^{\mu+1} v \rangle \nonumber \\
 &\ \  + \frac{\|B\|_{L^\infty([0,\infty) \times \Omega)}b_0^2 R^2}{\lambda K} \|\nabla \partial_t^\mu v \|^2_2 \nonumber + C \lambda Z_0(v)+ \langle \partial_t^\mu F_\lambda,  [h ; \nabla \partial_t^\mu v] \rangle \nonumber,
\end{align}
where $\sigma$ is the outward pointing unit normal vector of $\partial \Omega_\lambda$,
\begin{equation}\label{eq:h}
\phi(r)= \left\{
\begin{array}{ll}
b_0 , &\quad (r \leq \frac{R}{\lambda}) \\
\frac{b_0 R}{\lambda r} , &\quad ( r \geq \frac{R}{\lambda})
\end{array}
\right., \quad h(x)=x \phi(|x|)
\end{equation}
and $[h;\nabla g] :\mathbb{R}^d \to \mathbb{R}^d$ is defined by
\begin{equation}\nonumber
([h;\nabla g])^i(x) = h(x) \cdot \nabla g^{i}(x) \quad (i = 1,2, \cdots, d)
\end{equation}
for any $g:\mathbb{R}^d \to \mathbb{R}^d$.

\proof
Let $0 \leq \mu \leq L-1$ and $\lambda \leq 1$. First, we prove \eqref{eq:A}.
Applying $\partial_t^\mu$ to (DW) and taking inner product it by $\partial_t^{\mu+1} u$, we have
\begin{align}\label{eq:3.11}
 &\frac{d}{dt} E(\partial_t^\mu v) + \langle \partial_t^{\mu+1}  v, B_\lambda \partial_t^{\mu+1} v \rangle \\
 &= - \sum_{1 \leq \nu \leq \mu} \binom{\mu}{\nu} \langle \partial_t^{\mu+1} v, \partial_t^\nu B_\lambda \partial_t^{\mu-\nu +1} v \rangle + \langle \partial_t^{\mu+1} v, \partial_t^\mu F_\lambda \rangle \nonumber ,
\end{align}
where we use
\begin{equation}\label{eq:3.12}
\partial_t^\mu v = 0 \quad  on \quad \partial \Omega_\lambda.
\end{equation}
Note that even if $\mu =0$, \eqref{eq:3.11} holds in the sense of the first term in the right-hand side to be zero. 
Since $\lambda \leq 1$, the definition of $Z_0$ and condition $\bold{(B4)_\lambda}$ imply that
\begin{equation}\nonumber
- \sum_{1 \leq \nu \leq \mu} \binom{\mu}{\nu} \langle \partial_t^{\mu+1} v, \partial_t^\nu B_\lambda \partial_t^{\mu-\nu +1} v \rangle \leq C Z_0 \sum_{1 \leq \nu \leq \mu} \|\partial_t^\nu B_\lambda\|_\infty \leq C \lambda^2 Z_0,
\end{equation}
thus we obtain \eqref{eq:A}.

Second, we prove \eqref{eq:B}. Applying $\partial_t^\mu$ to (DW) and taking inner product it by $\partial_t^\mu u$, we have
\begin{align}
&\frac{d}{dt} \{ \langle \partial_t^\mu v, \partial_t^{\mu+1} v \rangle + \frac{1}{2} \langle \partial_t^\mu v, B_\lambda \partial_t^\mu v \rangle \} + \|\nabla \partial_t^\mu v \|^2_2 - \| \partial_t^{\mu+1} v\|^2_2 \nonumber \\
& =   \langle \partial_t^\mu v, \partial_t^{\mu+2} v \rangle + \langle \partial_t^\mu v, B_\lambda \partial_t^{\mu+1} v \rangle  + \|\nabla \partial_t^\mu v \|^2_2 +  \frac{1}{2} \langle \partial_t^\mu v, \partial_t B_\lambda \partial_t^\mu v \rangle  \nonumber \\
&= - \sum_{1 \leq \nu \leq \mu} \binom{\mu}{\nu} \langle \partial_t^\mu  v, \partial_t^\nu B_\lambda \partial_t^{\mu-\nu +1}  v \rangle + \frac{1}{2} \langle \partial_t^\mu v, \partial_t B_\lambda \partial_t^\mu v \rangle + \langle \partial_t^\mu  v, \partial_t^\mu F_\lambda \rangle \nonumber \\
&\leq C \lambda^2 Z_0  + \langle \partial_t^\mu  v, \partial_t^\mu F_\lambda \rangle \nonumber,
\end{align}
where we use assumption $\bold{(B2)_\lambda}$. It means that \eqref{eq:B} holds. 

Finally, we prove \eqref{eq:C}. Applying $\partial_t^\mu $ to (DW$)_\lambda$ and taking inner product (DW) by $[h ; \nabla \partial_t^\mu v]$ we obtain
\begin{align}\label{eq:3.13}
 &\langle \partial_t^{\mu+2}  v, [h ; \nabla \partial_t^\mu v] \rangle -\langle \triangle \partial_t^\mu v, [h ; \nabla \partial_t^\mu v] \rangle + \langle B_\lambda \partial_t^{\mu+1} v,[h ; \nabla \partial_t^\mu v] \rangle \\
 &= -\sum_{1 \leq \nu \leq \mu} \binom{\mu}{\nu} \langle \partial_t^\nu B_\lambda \partial_t^{\mu-\nu+1} v,[h ; \nabla \partial_t^\mu v] \rangle + \langle \partial_t^\mu F_\lambda, [h ; \nabla \partial_t^\mu v] \rangle \nonumber
\end{align}
Noting
\begin{equation}\label{eq:3.14}
\nabla \partial_t^\mu v^k = \sigma \cdot \nabla  \partial_t^\mu v^k \sigma \quad on \ \partial  \Omega_\lambda,
\end{equation}
we obtain
\begin{align*}
 & \langle \partial_t^{\mu+2}  v, [h ; \nabla \partial_t^\mu v] \rangle -\frac{d}{dt} \langle \partial_t^{\mu+1} v , [h ; \nabla \partial_t^\mu v] \rangle = - \sum_{i,k=1}^d \int_{\Omega_\lambda} \partial_t^{\mu+1} v^k \partial_i \partial_t^{\mu+1} v^k h^i dx \\
 &=- \frac{1}{2} \int_{\Omega_\lambda} h \cdot \nabla |\partial_t^{\mu+1}  v|^2 dx = \frac{1}{2} \int_{\Omega_\lambda} \textrm{div} h |\partial_t^{\mu+1} v|^2  dx
\end{align*}
and
\begin{eqnarray*}
 && - \left\langle \triangle \partial_t^\mu v, [h ; \nabla \partial_t^\mu v]  \right\rangle  \\
 &=& \sum_{k=1}^d \int_{\Omega_\lambda} \nabla \partial_t^\mu v^k \cdot \nabla ( h \cdot \nabla \partial_t^\mu v^k ) dx - \sum_{k=1}^d \int_{\partial \Omega_\lambda} \sigma \cdot \nabla \partial_t^\mu v^k h \cdot \nabla \partial_t^\mu v^k dS \\
 &=& \sum_{i,j,k=1}^d \int_{\Omega_\lambda} \partial_j \partial_t^\mu v^k  \partial_i \partial_t^\mu v^k \partial_ih^j  dx + \frac{1}{2} \int_{\Omega_\lambda} h \cdot \nabla  |\nabla \partial_t^\mu v|^2 dx \\
 && -  \int_{\partial \Omega_\lambda} h \cdot \sigma | \sigma \cdot \nabla \partial_t^\mu v|^2  dS \\
 &=& \sum_{i,j,k=1}^d \int_{\Omega_\lambda} \partial_j \partial_t^\mu v^k  \partial_i \partial_t^\mu v^k \partial_ih^j  dx - \frac{1}{2} \int_{\Omega_\lambda} \textrm{div} h |\nabla \partial_t^\mu v|^2 dx \\
 && -\frac{1}{2} \int_{\partial \Omega_\lambda} h \cdot \sigma | \sigma \cdot \nabla \partial_t^\mu v|^2  dS.
\end{eqnarray*}
Therefore we get from \eqref{eq:3.13} that
\begin{eqnarray}\label{eq:3.15}
 && \frac{d}{dt} \langle \partial_t^{\mu+1} v , [h ; \nabla \partial_t^\mu v] \rangle + \frac{1}{2} \int_{\Omega_\lambda} ( |\partial_t^{\mu+1} v|^2 - |\nabla \partial_t^\mu v|^2 ) \textrm{div} h  dx \\
 && + \sum_{i,j,k=1}^d \int_{\Omega_\lambda} \partial_j \partial_t^\mu v^k  \partial_i \partial_t^\mu v^k \partial_ih^j  dx  \nonumber \\
 &=& \frac{1}{2} \int_{\partial \Omega_\lambda} h \cdot \sigma | \sigma \cdot \nabla \partial_t^\mu v|^2  dS -  \langle B_\lambda \partial_t^{\mu+1} v,[h ; \nabla \partial_t^\mu v] \rangle \nonumber \\
 && - \sum_{1 \leq \nu \leq \mu} \binom{\mu}{\nu} \langle \partial_t^\nu B_\lambda \partial_t^{\mu-\nu+1} v,[h ; \nabla \partial_t^\mu v] \rangle + \langle \partial_t^\mu F_\lambda,  [h ; \nabla \partial_t^\mu v] \rangle. \nonumber
\end{eqnarray}
Now we remark that
\begin{equation}\label{eq:hes}
\left\{
\begin{array}{ll}
\displaystyle \frac{\partial h^{i}}{\partial x_j} = \delta_{ij} \phi(|x|) + \phi'(|x|) \frac{x_ix_j}{|x|}, &  \\
\displaystyle \textrm{div}h(x) = d \phi(|x|) + \phi'(|x|) |x|, & \\
\displaystyle \|h \|_\infty \leq  \frac{b_0 R}{\lambda}\quad and \quad \|\nabla h \|_\infty \leq 2 b_0.  &
\end{array}
\right.
\end{equation}
hold. Using \eqref{eq:hes} and $\phi' (r) \leq 0$, we obtain 
\begin{eqnarray*}
 && \sum_{i,j,k=1}^d \int_{\Omega_\lambda} \partial_j \partial_t^\mu v^k  \partial_i \partial_t^\mu v^k \partial_ih^j  dx\\
 &=& \sum_{i,k=1}^d \int_{\Omega_\lambda} |\partial_i \partial_t^\mu v^k|^2 \phi dx + \sum_{i,j,k=1}^d \int_{\Omega_\lambda} \partial_j \partial_t^\mu v^k  \partial_i \partial_t^\mu v^k \phi' \frac{x_i x_j}{|x|}  dx \\
 &=& \int_{\Omega_\lambda} |\nabla \partial_t^\mu v|^2 \phi dx + \sum_{k=1}^d \int_{\Omega_\lambda} |x \cdot \nabla \partial_t^\mu v^k |^2 \phi' \frac{1}{|x|}  dx \\
 &\geq&  \int_{\Omega_\lambda} \{ \phi + |x| \phi' \} |\nabla \partial_t^\mu v|^2 dx .
\end{eqnarray*}
These estimates and \eqref{eq:3.15} imply that
\begin{align}\label{eq:3.17}
& \frac{d}{dt} \langle \partial_t^{\mu+1} v, [h ; \nabla \partial_t^\mu v] \rangle  \\
& \quad + \int_{\Omega_\lambda}  \left( \frac{ (d \phi + |x| \phi')}{2} \right) |\partial_t^{\mu+1} v|^2 dx + \int_{\Omega_\lambda}  \left( \frac{(2-d)\phi + |x| \phi'}{2}  \right) |\nabla \partial_t^\mu v |^2 dx \nonumber  \\
&\leq \frac{1}{2} \int_{\partial \Omega_\lambda} h \cdot \sigma | \sigma \cdot \nabla \partial_t^\mu v|^2  dS - \langle B_\lambda \partial_t^{\mu+1} v, [h ; \nabla \partial_t^\mu v] \rangle  \nonumber \\
& \quad - \sum_{1 \leq \nu \leq \mu} \binom{\mu}{\nu} \langle\partial_t^\nu B_\lambda \partial_t^{\mu-\nu+1} v, [h ; \nabla \partial_t^\mu v] \rangle  + \langle \partial_t^\mu F_\lambda,  [h ; \nabla \partial_t^\mu v] \rangle. \nonumber
\end{align}
Let we estimate for the right side of \eqref{eq:3.17}. We calculate
\begin{align}\label{eq:3.18}
|\langle B_\lambda \partial_t^{\mu+1} v, [h; \nabla \partial_t^\mu v \rangle| &= |\langle \sqrt{B_\lambda} \partial_t^{\mu+1} v, \sqrt{B_\lambda}[ h ; \nabla \partial_t^\mu v] \rangle| \\
\nonumber &\leq  \frac{K}{4}\| \sqrt{B_\lambda}\partial_t^{\mu+1} v\|_2^2 + \frac{1}{K}\|\sqrt{B_\lambda} [h; \nabla \partial_t^\mu v] \|_2^2 \\
\nonumber &\leq \frac{K}{4} \langle \partial_t^{\mu+1} v, B_\lambda\partial_t^{\mu+1} v \rangle + \frac{\|B\|_{L^\infty([0,\infty) \times \Omega)} b_0^2 R^2}{\lambda K} \|\nabla \partial_t^\mu v \|^2_2 ,
\end{align}
\begin{eqnarray}\label{eq:3.18}
&& \left| \sum_{1 \leq \nu \leq \mu} \binom{\mu}{\nu} \langle \partial_t^\nu B_\lambda \partial_t^{\mu-\nu+1} v, [h ; \nabla \partial_t^\mu v ] \rangle \right| \\
\nonumber &\leq& C  \sum_{1 \leq \nu \leq \mu}  \| \partial_t^\nu  B_\lambda \|_\infty \| h \|_\infty \|\partial_t^{\mu-\nu+1} v \|_2 \|\nabla \partial_t^\mu v \|_2  \leq C \lambda Z_0.
\end{eqnarray}
where we use $\bold{(B1)_\lambda}, \bold{(B4)_\lambda}$ and \eqref{eq:hes}.
Combining estimates \eqref{eq:3.17} - \eqref{eq:3.18}, we get \eqref{eq:C}.
This completes the proof of Lemma \ref{lem:3.4}.

\qed
\end{lem}

We define $G(v(t))$ by
\begin{align}\label{eq:Gmu}
G(v(t)) &= \frac{C_0}{2\lambda} Z_0(v(t)) + \frac{b_0 (2d-1)}{4} \sum_{\mu=0}^{L-1} \langle \partial_t^{\mu+1} v(t), \partial_t^\mu v(t) \rangle  \\
\nonumber & + \frac{b_0(2d-1)}{8} \sum_{\mu=0}^{L-1} \langle \partial_t^\mu v(t), B_\lambda(t) \partial_t^\mu v(t) \rangle  + \sum_{\mu=0}^{L-1} \langle \partial_t^{\mu+1} v(t), [h ; \nabla \partial_t^\mu v(t)] \rangle,
\end{align}
where
\begin{equation}\label{eq:C0k}
C_0 = \max\left\{4b_0R + \frac{C_1b_0(2d-1)}{2} , \ d , \ \|B\|_{L^\infty([0,\infty) \times \Omega)} b_0^2 R^2 \times \frac{8}{b_0}   \right\}.
\end{equation}

\begin{lem}\label{lem:3.5}
There exists a constant $0< \lambda <1$ such that the local solution $v$ to {\rm(DW$)_\lambda$} satisfies
\begin{align}\label{eq:3.22}
&\frac{d}{dt} G(v(t)) + \frac{b_0}{16} Z_0(v(t)) \\
 &\leq \frac{1}{2} \sum_{\mu=0}^{L-1} \int_{\partial \Omega_\lambda} h \cdot \sigma | \sigma \cdot \nabla \partial_t^\mu v|^2  dS + \bar{C}_\lambda \sum_{\mu=0}^{L-1}  \langle \partial_t^\mu F_\lambda, \partial_t^{\mu+1} v + \partial_t^\mu v + [h ; \nabla \partial_t^\mu v] \rangle \nonumber,
\end{align}
where $\displaystyle \bar{C}_\lambda = \frac{C_0}{\lambda} + \frac{b_0(2d-1)}{4} + 1$.

\proof 
Let $\mu \leq L-1$ and $K \geq \frac{d}{\lambda}$. Calculating $K \times \eqref{eq:A} + \frac{b_0(2d-1)}{4} \times \eqref{eq:B} + \eqref{eq:C}$,  we have
\begin{align}\label{eq:3.23}
 & \frac{d}{dt} \left\{ K E(\partial_t^\mu v) + \frac{b_0 (2d-1)}{4} \langle \partial_t^{\mu+1}  v, \partial_t^\mu v \rangle \right. \\
 & \hspace{30mm} + \left. \frac{b_0(2d-1)}{8}\langle \partial_t^\mu v, B_\lambda \partial_t^\mu v \rangle  + \langle \partial_t^{\mu+1} v, [h; \nabla \partial_t^\mu v] \rangle \right\} \nonumber \\
 &+ K \langle \partial_t^{\mu+1} v, B_\lambda \partial_t^{\mu+1} v \rangle - \frac{b_0(2d-1)}{4} \|\partial_t^{\mu+1} v\|^2_2 + \int_{\Omega_\lambda} \left( \frac{d \phi + |x| \phi'}{2} \right) |\partial_t^{\mu+1} v|^2 dx \nonumber \\
 &+ \frac{b_0(2d-1)}{4} \|\nabla \partial_t^\mu v \|^2_2 + \int_{\Omega_\lambda} \left( \frac{(2-d) \phi + |x| \phi'}{2} \right) |\nabla \partial_t^\mu v|^2  dx \nonumber \\
 &\leq  \frac{K}{4} \langle \partial_t^{\mu+1} v, B_\lambda \partial_t^{\mu+1} v \rangle + \frac{\|B\|_{L^\infty([0,\infty) \times \Omega)}b_0^2 R^2}{\lambda K} \|\nabla \partial_t^\mu v \|^2_2 \nonumber \\
 & + \frac{1}{2} \int_{\partial \Omega_\lambda} h \cdot \sigma | \sigma \cdot \nabla \partial_t^\mu v|^2  dS  \nonumber + C\lambda^2 Z_0 \left( K + \frac{b_0(2d-1)}{4}  + \frac{1}{\lambda} \right)  \nonumber  \\
 &+ \left( K + \frac{b_0(2d-1)}{4}  + 1 \right) \langle \partial_t^\mu F_\lambda, \partial_t^{\mu+1} v + \partial_t^\mu v + [h ; \nabla \partial_t^\mu v] \rangle . \nonumber
\end{align}
Using \textbf{(B3$)_\lambda$}, $K \geq \frac{d}{\lambda}$, $\phi \geq 0$ and 
\begin{equation}\nonumber
r \phi'(r) =\left\{
\begin{array}{ll}
0 , &\quad (r < \frac{R}{\lambda}) \\
-\phi(r) , &\quad ( r > \frac{R}{\lambda}),
\end{array}
\right.
\end{equation}
we obtain
\begin{eqnarray}\label{eq:3.24}
&& K \langle \partial_t^{\mu+1} v,B_\lambda \partial_t^{\mu+1} v \rangle -\frac{b_0(2d-1)}{4}\| \partial_t^{\mu+1} v \|^2_2 + \int_{\Omega_\lambda} \left\{ \frac{d\phi+|x| \phi'}{2} \right\} |\partial_t^{\mu+1} v|^2 dx \nonumber \\
\nonumber &\geq& \frac{K}{2} \langle \partial_t^{\mu+1} v,B_\lambda\partial_t^{\mu+1} v \rangle + \int_{U_{\frac{R}{\lambda}}} \left( -\frac{b_0(2d-1)}{4} + \frac{db_0}{2}\right) |\partial_t^{\mu+1} v|^2 dx \\
\nonumber && + \int_{|x| \geq \frac{R}{\lambda}} \left(\frac{\lambda b_0 K}{2} - \frac{b_0(2d-1)}{4} + \frac{d-1}{2} \phi  \right) |\partial_t^{\mu+1} v|^2  dx \\
&\geq& \frac{K}{2} \langle \partial_t^{\mu+1} v,B_\lambda\partial_t^{\mu+1} v \rangle +  \frac{b_0}{4} \|\partial_t^{\mu+1} v \|^2_2
\end{eqnarray}
and
\begin{eqnarray}\label{eq:3.25}
&& \frac{b_0(2d-1)}{4} \|\nabla \partial_t^\mu v \|_2^2 + \int_{\Omega_\lambda} \frac{(2-d) \phi + |x| \phi'}{2} |\nabla \partial_t^\mu v |^2 dx \\
\nonumber &=& \int_{U_{\frac{R}{\lambda}}} \left( \frac{b_0(2d-1)}{4} + \frac{(2-d) b_0}{2}  \right) |\nabla \partial_t^\mu v |^2 dx \\
\nonumber && + \int_{|x| \geq\frac{R}{\lambda}} \left(\frac{b_0(2d-1)}{4} + \frac{(1-d)}{2} \frac{b_0 R}{\lambda |x|}  \right) |\nabla \partial_t^\mu v |^2dx \geq \frac{b_0}{4} \| \nabla \partial_t^\mu v \|^2_2.
\end{eqnarray}
Combining \eqref{eq:3.23}, \eqref{eq:3.24} and \eqref{eq:3.25}, we have
\begin{align*}
 & \frac{d}{dt} \left\{ K E(\partial_t^\mu v) + \frac{b_0 (2d-1)}{4} \langle \partial_t^{\mu+1}  v, \partial_t^\mu v \rangle \right. \\
 & \hspace{30mm} + \left. \frac{b_0(2d-1)}{8}\langle \partial_t^\mu v, B_\lambda \partial_t^\mu v \rangle  + \langle \partial_t^{\mu+1} v, [h; \nabla \partial_t^\mu v] \rangle \right\} \nonumber \\
 &+ \frac{b_0}{4} \{ \|\partial_t^{\mu+1} v \|^2_2 + \|\nabla \partial_t^\mu v\|^2_2 \} \nonumber \\
 & \leq \frac{\|B\|_{L^\infty([0,\infty) \times \Omega)}b_0^2 R^2}{\lambda K} \|\nabla \partial_t^\mu v \|^2_2 + \frac{1}{2} \int_{\partial \Omega_\lambda} h \cdot \sigma | \sigma \cdot \nabla \partial_t^\mu v|^2  dS \nonumber \\
 &  + C \lambda^2 Z_0(v) \left(K + \frac{b_0(2d-1)}{4}  + \frac{1}{\lambda} \right)  \nonumber  \\
 & +\left( K + \frac{b_0(2d-1)}{4}  + 1\right) \langle \partial_t^\mu F_\lambda, \partial_t^{\mu+1} v + \partial_t^\mu v + [h ; \nabla \partial_t^\mu v] \rangle. \nonumber
\end{align*}
Let $K= \frac{C_0}{\lambda}$ and sum up $\mu$ from $0$ to $L-1$. Then we get
\begin{align*}
&\frac{d}{dt} G(v) + \frac{b_0}{8} Z_0(v) \leq C \lambda Z_0 (v)  \\
& + \frac{1}{2} \sum_{\mu=0}^{L-1} \int_{\partial \Omega_\lambda} h \cdot \sigma | \sigma \cdot \nabla \partial_t^\mu v|^2  dS + \bar{C}_\lambda \sum_{\mu=0}^{L-1}  \langle \partial_t^\mu F_\lambda, \partial_t^{\mu+1} v + \partial_t^\mu v + [h ; \nabla \partial_t^\mu v] \rangle, \nonumber
\end{align*}
where we use \eqref{eq:C0k}. We choose $\lambda$ satisfies small such that $\lambda$ satisfies $C \lambda \leq \frac{b_0}{16}$. Then we obtain
\begin{align*}
&\frac{d}{dt} G (v) + \frac{b_0}{16} Z_0(v) \\
& \leq \frac{1}{2} \sum_{\mu=0}^{L-1} \int_{\partial \Omega} h \cdot \sigma | \sigma \cdot \nabla \partial_t^\mu v|^2  dS + \bar{C}_\lambda \sum_{\mu=0}^{L-1} \langle \partial_t^\mu F_\lambda, \partial_t^{\mu+1} v + \partial_t^\mu v + [h ; \nabla \partial_t^\mu v] \rangle. \nonumber 
\end{align*}
This completes the proof of Lemma \ref{lem:3.5}.
\qed
\end{lem}

We choose the $\lambda$ sufficiently small to hold the Lemma \ref{lem:3.5}. Next Lemma are the estimates of the nonlinear terms.

\begin{lem}\label{lem:3.6}
Let  $\mu \leq L-1$. Then there exists a constant $C_\lambda>0$ such that the local solution $v \in X_\delta^T$ to {\rm(DW$)_\lambda$} satisfies
\begin{align}\label{eq:A'}
 &\langle \partial_t^{\mu+1} v, \partial_t^\mu F_\lambda  \rangle \leq C_\lambda \delta Z_0 - \frac{1}{2} \sum_{i,j=1}^d \sum_{1 \leq a,b \leq d} \frac{d}{dt} \int_{\Omega_\lambda} \partial_t^\mu \partial_a v^i \partial_t^\mu \partial_b v^j c^{ab}_{ij} (\partial v) dx \\
 & \qquad \qquad + \sum_{i,j=1}^d \frac{1}{2} \frac{d}{dt} \int_{\Omega_\lambda} \partial_t^{\mu+1} v^i \partial_t^{\mu+1} v^j c^{00}_{ij} (\partial v) dx, \nonumber
\end{align}
\begin{align}\label{eq:B'}
\langle \partial_t^{\mu}  v, \partial_t^\mu F_\lambda  \rangle \leq C_\lambda \delta Z_0 + \sum_{i,j=1}^d \sum_{0 \leq b \leq d} \frac{d}{dt} \int_{\Omega_\lambda} \partial_t^\mu v^i \partial_t^\mu \partial_b v^j c^{0b}_{ij} (\partial v) dx
\end{align}
and
\begin{align}\label{eq:C'}
 &\langle [h; \nabla \partial_t^\mu v], \partial_t^\mu F_\lambda  \rangle \leq C_\lambda \delta Z_0 + C_\lambda \delta \int_{\partial \Omega_\lambda }|h \cdot \sigma| |\sigma \cdot \nabla \partial_t^\mu v |^2 dS \\
 &\qquad \qquad + \sum_{i,j=1}^d \sum_{0 \leq b \leq d} \frac{d}{dt} \int_{\Omega_\lambda} h \cdot \nabla \partial_t^\mu v^i \partial_t^\mu \partial_b v^j c^{0b}_{ij} (\partial v)  dx \nonumber.
\end{align}
\proof
We may assume that $\delta<1$. First, we prove \eqref{eq:A'}. Using Lemma \ref{lem:3.3} and Lemma \ref{lem:A1}, we have
\begin{equation}\nonumber
\langle \partial_t^{\mu+1} v, \partial_t^\mu \tilde{F}_\lambda  \rangle \leq \|\partial_t^{\mu+1} v\|_2 \|\partial_t^\mu \tilde{F}_\lambda \|_2 \leq C_\lambda \delta Z \leq C_\lambda \delta Z_0.
\end{equation}
Furthermore we calculate that
\begin{align*}
 & \sum_{i,j=1}^d \sum_{0 \leq a,b \leq d} \int_{\Omega_\lambda} \partial_t^{\mu+1}v^i \partial_t^\mu(c^{ab}_{ij}(\partial v) \partial_a \partial_b v^j) dx \\
 &= \sum_{0 \leq \nu \leq \mu-1}\binom{\mu}{\nu} \sum_{i,j=1}^d \sum_{0 \leq a,b \leq d} \int_{\Omega_\lambda} \partial_t^{\mu+1} v^i \partial_t^\nu \partial_a \partial_b v^j \partial_t^{\mu-\nu}(c^{ab}_{ij} (\partial v)) dx  \\
 & - \sum_{i,j=1}^d \sum_{0 \leq a,b \leq d} \int_{\Omega_\lambda} \partial_t^{\mu+1} v^i \partial_t^\mu \partial_b v^j \partial_a( c^{ab}_{ij}(\partial v) ) dx\\
 & - \sum_{i,j=1}^d \sum_{0 \leq a,b \leq d} \int_{\Omega_\lambda} \partial_t^{\mu+1} \partial_a v^i \partial_t^\mu \partial_b v^j c^{ab}_{ij}(\partial v) dx\\
 & + \sum_{i,j=1}^d \sum_{0 \leq a,b \leq d} \int_{\Omega_\lambda} \partial_a ( \partial_t^{\mu+1} v^i \partial_t^\mu \partial_b v^j c^{ab}_{ij} (\partial v) ) dx \\
 &= J_1 + J_2  + \sum_{i,j=1}^d \sum_{0 \leq a,b \leq d} \frac{1}{2} \int_{\Omega_\lambda} \partial_t^\mu \partial_a v^i \partial_t^\mu \partial_b v^j \partial_t( c^{ab}_{ij} (\partial v)) dx \\
 & - \sum_{i,j=1}^d \sum_{0 \leq a,b \leq d} \frac{1}{2} \frac{d}{dt} \int_{\Omega_\lambda} \partial_t^\mu \partial_a v^i \partial_t^\mu \partial_b v^j c^{ab}_{ij} (\partial v) dx  \\
 & +  \sum_{i,j=1}^d \sum_{0 \leq b \leq d}  \frac{d}{dt} \int_{\Omega_\lambda} \partial_t^{\mu+1} v^i \partial_t^\mu \partial_b v^j c^{0b}_{ij}(\partial v) dx\\
 & = J_1 + J_2 + J_3 - \sum_{i,j=1}^d \sum_{1 \leq a,b \leq d} \frac{1}{2} \frac{d}{dt} \int_{\Omega_\lambda} \partial_t^\mu \partial_a v^i \partial_t^\mu \partial_b v^j c^{ab}_{ij} (\partial v) dx\\
 & \qquad \qquad + \sum_{i,j=1}^d \frac{1}{2} \frac{d}{dt} \int_{\Omega_\lambda} \partial_t^{\mu+1} v^i \partial_t^{\mu+1} v^j c^{00}_{ij} (\partial v) dx .
\end{align*}
We can estimate $|J_k| \leq C_\lambda \delta Z_0 \ (k=1,2,3)$ from Lemma \ref{lem:3.3} and Lemma \ref{lem:A2}. Therefore we get \eqref{eq:A'}.

Second, we prove \eqref{eq:B'}. Using Lemma \ref{lem:3.3} and Lemma \ref{lem:A1}, we have
\begin{equation}\nonumber
\langle \partial_t^{\mu} v, \partial_t^\mu \tilde{F}_\lambda  \rangle \leq \|\partial_t^{\mu} v\|_2 \|\partial_t^\mu \tilde{F}_\lambda \|_2 \leq C_\lambda \delta Z \leq C_\lambda \delta Z_0.
\end{equation}
Furthermore we calculate that
\begin{align*}
 & \sum_{i,j=1}^d \sum_{0 \leq a,b \leq d} \int_{\Omega_\lambda} \partial_t^{\mu} v^i \partial_t^\mu(c^{ab}_{ij}(\partial v) \partial_a \partial_b v^j) dx \\
 &= \sum_{0 \leq \nu \leq \mu-1} \binom{\mu}{\nu} \sum_{i,j=1}^d \sum_{0 \leq a,b \leq d} \int_{\Omega_\lambda} \partial_t^{\mu} v^i \partial_t^\nu \partial_a \partial_b v^j \partial_t^{\mu-\nu}(c^{ab}_{ij} (\partial v)) dx  \\
 & - \sum_{i,j=1}^d \sum_{0 \leq a,b \leq d} \int_{\Omega_\lambda} \partial_a \partial_t^{\mu} v^i \partial_t^\mu \partial_b v^j c^{ab}_{ij} (\partial v) dx \\
 & - \sum_{i,j=1}^d \sum_{0 \leq a,b \leq d}  \int_{\Omega_\lambda} \partial_t^\mu v^i \partial_t^\mu \partial_b v^j \partial_a (c^{ab}_{ij}(\partial v)) dx\\
 & + \sum_{i,j=1}^d \sum_{0 \leq b \leq d} \frac{d}{dt} \int_{\Omega_\lambda} \partial_t^\mu v^i \partial_t^\mu \partial_b v^j c^{0b}_{ij} (\partial v) dx \\
 & = J_4 + J_5 + J_6 + \sum_{i,j=1}^d \sum_{0\leq b \leq d} \frac{d}{dt} \int_{\Omega_\lambda} \partial_t^\mu v^i \partial_t^\mu \partial_b v^j c^{0b}_{ij} (\partial v) dx.
\end{align*}
We can estimate $|J_k| \leq C_\lambda \delta Z_0 \ (k=4,5,6)$ from Lemma \ref{lem:3.3} and Lemma \ref{lem:A2}. Therefore we get \eqref{eq:B'}.

Finally, we prove \eqref{eq:C'}, Using Lemma \ref{lem:3.3} and Lemma \ref{lem:A1}, we have
\begin{equation}\nonumber
\langle [h; \nabla \partial_t^\mu v], \partial_t^\mu \tilde{F}_\lambda  \rangle \leq \|h\|_\infty \|\nabla \partial_t^{\mu} v\|_2 \|\partial_t^\mu \tilde{F}_\lambda \|_2 \leq C_\lambda \delta Z \leq C_\lambda \delta Z_0.
\end{equation}
We calculate that
\begin{align*}
 & \sum_{i,j=1}^d \sum_{0 \leq a,b \leq d} \int_{\Omega_\lambda} h \cdot \nabla \partial_t^\mu v^i \partial_t^\mu(c^{ab}_{ij}(\partial v) \partial_a \partial_b v^j) dx \\
 &= \sum_{0 \leq \nu \leq \mu-1}\sum_{i,j=1}^d \sum_{0 \leq a,b \leq d} \int_{\Omega_\lambda} h \cdot \nabla \partial_t^\mu v^i \partial_t^\nu \partial_a \partial_b v^j \partial_t^{\mu-\nu}(c^{ab}_{ij} (\partial v)) dx  \\
 &- \sum_{i,j=1}^d \sum_{0 \leq a,b \leq d} \int_{\Omega_\lambda}  h \cdot \nabla \partial_t^\mu v^i \partial_t^\mu \partial_b v^j \partial_a ( c^{ab}_{ij} (\partial v) )dx \\
 &- \sum_{i,j=1}^d \sum_{0 \leq a,b \leq d} \int_{\Omega_\lambda} \partial_a h \cdot \nabla \partial_t^\mu v^i \partial_t^\mu \partial_b v^j c^{ab}_{ij} (\partial v) dx \\
 &+ \sum_{i,j=1}^d \sum_{0 \leq a,b \leq d} \frac{1}{2} \int_{\Omega_\lambda}  h \cdot \nabla( c^{ab}_{ij} (\partial v) ) \partial_t^\mu \partial_a v^i \partial_t^\mu \partial_b v^jdx \\
 &+ \sum_{i,j=1}^d \sum_{0 \leq a,b \leq d} \frac{1}{2} \int_{\Omega_\lambda} ({\rm div} h) \partial_t^\mu \partial_a v^i \partial_t^\mu \partial_b v^j c^{ab}_{ij} (\partial v) dx\\
 &+ \sum_{i,j=1}^d \sum_{1 \leq a,b \leq d} \int_{\partial \Omega_{\lambda}} \sigma_a h \cdot \nabla \partial_t^\mu v^i \partial_t^\mu \partial_b v^j c^{ab}_{ij} (\partial v)  dS\\
 &- \sum_{i,j=1}^d \sum_{1 \leq a,b \leq d} \frac{1}{2} \int_{\partial \Omega_\lambda} h \cdot \sigma \partial_t^\mu \partial_a v^i \partial_t^\mu \partial_b v^j c^{ab}_{ij} (\partial v) dS \\
 &+ \sum_{i,j=1}^d \sum_{0 \leq b\leq d} \frac{d}{dt} \int_{\Omega_\lambda} h \cdot \nabla \partial_t^\mu v^i \partial_t^\mu \partial_b v^j c^{0b}_{ij} (\partial v) dx  \\
 & = J_7 + J_8 + J_9 + J_{10} + J_{11} +J_{12} + J_{13} \\
 & \qquad + \sum_{i,j=1}^d \sum_{0 \leq b \leq d} \frac{d}{dt} \int_{\Omega_\lambda} h \cdot \nabla \partial_t^\mu v^i \partial_t^\mu \partial_b v^j c^{0b}_{ij} (\partial v)  dx.
\end{align*}
We can estimate $|J_k| \leq C_\lambda \delta Z_0 \ (k=7,8,9,10,11)$ from Lemma \ref{lem:3.3} and Lemma \ref{lem:A2}.
Moreover using \eqref{eq:3.14} and
\begin{equation}\nonumber
\|\partial v \|_{L^\infty(\partial \Omega_\lambda)} \leq C_\lambda \|\partial v \|_{{H^{\left[\frac{d-1}{2}\right]+1}}(\partial \Omega_\lambda)} \leq C_\lambda \|\partial v \|_{{H^{\left[\frac{d-1}{2}\right]+2}}} \leq C_\lambda Z(v) 
\end{equation}
(The second inequality is the trace theorem, see for instance:\cite{Lion}), we get 
\begin{align*}
 J_{k} &\leq C \sum_{i,j=1}^d \sum_{1 \leq a,b \leq d} \|c^{ab}_{ij}(\partial v) \|_{L^\infty(\partial \Omega_\lambda)} \int_{\partial \Omega_{\lambda}} |h\cdot \sigma | |\sigma \cdot \nabla \partial_t^\mu v|^2 dS \\
 &\leq C_\lambda \delta \int_{\partial \Omega_{\lambda}} |h\cdot \sigma | |\sigma \cdot \nabla \partial_t^\mu v|^2 dS \qquad (k=12,13).
\end{align*}
Therefore we get \eqref{eq:C'}. This completes the proof of Lemma \ref{lem:3.6}.
\qed
\end{lem}

We define $\tilde{G}$ as follows:
\begin{align*}
 \tilde{G}(v(t))=& G(v(t)) + \bar{C}_\lambda \sum_{\mu=0}^{L-1} \sum_{i,j=1}^d \left\{ \frac{1}{2} \sum_{1 \leq a,b \leq d} \int_{\Omega_\lambda} \partial_t^\mu \partial_a v^i \partial_t^\mu \partial_b v^j c^{ab}_{ij} (\partial v) dx  \right.\\
 & - \frac{1}{2} \int_{\Omega_\lambda} \partial_t^{\mu+1} v^i \partial_t^{\mu+1} v^j c^{00}_{ij} (\partial v) dx \nonumber \\
 & \left. - \sum_{0 \leq b \leq d} \left\{ \int_{\Omega_\lambda} \partial_t^\mu v^i \partial_t^\mu \partial_b v^j c^{0b}_{ij} (\partial v) dx + \int_{\Omega_\lambda} h \cdot \nabla \partial_t^\mu v^i \partial_t^\mu \partial_b v^j c^{0b}_{ij} (\partial v) dx \right\} \right\}.
\end{align*}
Then the following lemma holds.
\begin{lem}\label{lem:3.7}
Let $v \in X_\delta^T$ be the solution to {\rm(DW$)_\lambda$} and $\lambda$ is sufficiently small to hold the Lemma \ref{lem:3.5}.
Then there exist a $\delta=\delta(\lambda)$ such that v satisfy 
\begin{equation}\label{eq:3.29}
\frac{d}{dt} \tilde{G}(v(t)) + \frac{b_0}{32} Z_0(v(t)) \leq 0 \qquad (t \in [0,T])
\end{equation}
and
\begin{equation}\label{eq:3.30}
\tilde{G} (v(t)) \cong_\lambda \|v(t)\|^2_2 + Z_0 (v(t)) \qquad (t \in [0,T]),
\end{equation}
where comparability constant is independent of $\delta,t$ and $T$.
\proof
First, we prove \eqref{eq:3.29}. Using Lemma \ref{lem:3.5} and Lemma \ref{lem:3.6}, we obtain
\begin{align*}
&\frac{d}{dt}\tilde{G}(v(t)) + \frac{b_0}{16} Z_0 \\
&\leq C_\lambda \delta Z_0 + \frac{1}{2} \sum_{\mu=0}^{L-1}\int_{\partial \Omega_\lambda} h \cdot \sigma |\sigma \cdot \nabla \partial_t^\mu v|^2 dS + C_\lambda \delta \sum_{\mu=0}^{L-1} \int_{\partial \Omega_\lambda} |h \cdot \sigma| |\sigma \cdot \nabla \partial_t^\mu v|^2 dS.
\end{align*}
Because of $\mathbb{R}^d/ \Omega_\lambda$ is star shaped, it holds that $h\cdot \sigma \leq 0$ on $\partial \Omega_\lambda$. Then we can choose $\delta$ sufficiently small depend on $\lambda$ such that \eqref{eq:3.29} holds.

Next, we prove \eqref{eq:3.30}. It follows from \eqref{eq:2.1} that 
\begin{eqnarray}\label{eq:3.31}
&& \left| \frac{b_0(2d-1)}{4}\sum_{\mu=0}^{L-1} \langle \partial_t^\mu v , \partial_t^{\mu+1} v \rangle \right| \\
\nonumber &\leq& \frac{b_0(2d-1)}{4} \sum_{\mu=0}^{L-1} \left\{ \frac{ \lambda}{4C_1}\| \partial_t^\mu v \|^2_2 + \frac{C_1}{\lambda} \|\partial_t^{\mu+1}  v \|^2_2  \right\} \\
\nonumber &\leq& \frac{b_0(2d-1)}{4} \sum_{\mu=0}^{L-1} \left\{ \frac{1}{4} \langle \partial_t^\mu v ,  B_\lambda \partial_t^\mu v \rangle + \frac{1}{4 \lambda}  \|\nabla \partial_t^\mu v \|^2_2 + \frac{C_1}{\lambda}\|\partial_t^{\mu+1} v\|^2_2 \right\} \\
\nonumber &\leq& \frac{b_0(2d-1)}{16} \sum_{\mu=0}^{L-1} \langle \partial_t^\mu v , B_\lambda \partial_t^\mu v \rangle + \frac{C_1b_0(2d-1)}{4\lambda} Z_0
\end{eqnarray}
and
\begin{eqnarray}\label{eq:3.32}
\left| \sum_{\mu=0}^{L-1} \langle \partial_t^{\mu+1} v, [h ;\nabla \partial_t^\mu v] \rangle \right| \leq \sum_{\mu=0}^{L-1} \|\partial_t^{\mu+1} v \|_2 \|\nabla \partial_t^\mu v\|_2 \|h\|_\infty \leq \frac{b_0R}{\lambda} Z_0.
\end{eqnarray}
Using (\ref{eq:3.31}), (\ref{eq:3.32}), \eqref{eq:C0k} and Lemma \ref{lem:Poincare}, we have
\begin{eqnarray}\label{eq:3.33}
&&G(v) \\
&\geq& \nonumber \frac{1}{\lambda} \left( \frac{C_0}{2} - b_0R - \frac{C_1b_0(2d-1)}{4}  \right) Z_0 + \frac{b_0(2d-1)}{16} \sum_{\mu=0}^{L-1} \langle \partial_t^\mu v,B_\lambda \partial_t^\mu v \rangle\\
&\geq& \frac{b_0R}{\lambda} Z_0 + \frac{b_0(2d-1)}{16} \sum_{\mu=0}^{L-1} \langle \partial_t^\mu v,B_\lambda \partial_t^\mu v \rangle \geq C_\lambda (\|v\|^2_2 + Z_0)  \nonumber.
\end{eqnarray}
On the other hand, Lemma \ref{prop:sobolev}, Lemma \ref{lem:3.3} and Lemma \ref{lem:A2} imply that 
\begin{equation}\nonumber 
\left| \int_{\Omega_\lambda} \partial_t^\mu \partial_a v^i \partial_t^\mu \partial_b v^j c^{ab}_{ij} (\partial v) dx \right| \leq \|\partial_a \partial_t^\mu v^i \|_2 \| \partial_t^\mu \partial_b v^j c^{ab}_{ij}(\partial v)\|_2 \leq C_\lambda \delta Z_0,
\end{equation}
\begin{equation}\nonumber 
\left| \int_{\Omega_\lambda} \partial_t^{\mu+1} v^i \partial_t^{\mu+1} v^j c^{00}_{ij} (\partial v) dx \right| \leq \|\partial^{\mu+1} v^i \|_2 \|\partial_t^{\mu+1} v^j c^{00}_{ij}(\partial v) \|_2 \leq C_\lambda \delta Z_0, \nonumber \\
\end{equation}
\begin{align*}
& \left| \int_{\Omega_\lambda} \partial_t^\mu v^i \partial_t^\mu \partial_b v^j c^{0b}_{ij} (\partial v) dx \right| \leq \|\partial_t^\mu v^i\|_2 \|\partial_t^\mu \partial_b v^j c^{0b}_{ij} (\partial v)\|_2 \leq C_\lambda \delta Z_0
\end{align*}
and
\begin{equation}\nonumber
\left| \int_{\Omega_\lambda} h \cdot \nabla \partial_t^\mu v^i \partial_t^\mu \partial_b v^j c^{0b}_{ij} (\partial v) dx \right| \leq \|h\|_\infty \|\nabla \partial_t^\mu v^i \|_2  \|\partial_t^\mu \partial_b v^jc^{0b}_{ij} (\partial v)\|_2 \leq C_\lambda \delta Z_0.
\end{equation}
Since these estimates and \eqref{eq:3.33} imply that we can choose $\delta = \delta (\lambda)$ to hold $\tilde{G}(v) \geq C_\lambda (\|v\|_2^2+Z_0)$.
It is clear that $\tilde{G}(v(t)) \leq C_\lambda' (\|v\|_2^2+Z_0)$ is true. Thus it holds that \eqref{eq:3.30}. 
This completes the proof of Lemma \ref{lem:3.7}.
\qed
\end{lem}

\subsection*{Proof of Proposition \ref{prop:3.2}}
Let $\lambda$ and $\delta$ be sufficiently small to hold Lemma \ref{lem:3.5} and Lemma \ref{lem:3.7}.
Integrating \eqref{eq:3.29} over $[0,t]$, we have
\begin{equation}\nonumber
\tilde{G}(v(t)) + \frac{b_0}{32} \int_0^t Z_0(v(s))ds \leq \tilde{G}(v(t))|_{t=0}.
\end{equation}
Since \eqref{eq:3.30} imply that 
\begin{align*}
&\quad \|v(t)\|^2_2 + Z_0(v(t)) + \int_0^t Z_0(v(s)) ds \leq C_\lambda(\|v_0\|^2_2 + Z_0(v(t))|_{t=0}),
\end{align*}
furthermore using Lemma \ref{lem:3.3}, we get \eqref{eq:prop3.2}.
This completes the proof of Proposition \ref{prop:3.2}.
\qed

\section{Decay Estimate}

In this section, we prove Theorem \ref{thm:main2}. 
In what follows, $\lambda$ and $\delta$ be sufficiently small to hold Theorem \ref{thm:main}.
Let $v \in X_{\delta}$ be the solution to (DW$)_\lambda$.
Since \eqref{eq:3.29} implies 
\begin{equation}\nonumber
\frac{d}{dt} \{ (1+t) \tilde{G}(v(t)) \} = \tilde{G}(v(t)) + (1+t) \frac{d}{dt} \tilde{G}(v(t)) \leq \tilde{G}(v(t)) - \frac{b_0}{32} (1+t) Z_0 (v(t)).
\end{equation}
Integrating the above estimate over $[0,t]$ and using \eqref{eq:3.30}, Lemma \ref{lem:Poincare} and Proposition \ref{prop:3.2}, we obtain 
\begin{align}\label{eq:4.1}
 &(1+t)\{ \|v(t)\|^2_2 + Z(v(t)) \} + \int_0^t (1+s) Z(v(s)) ds \\
 &\quad \leq C_\lambda \|(v_0,v_1)\|_{H^L \times H^{L-1}}^2 + C_\lambda \int_0^t \langle v(s) , B_\lambda(s) v(s)  \rangle ds \nonumber.
\end{align}
We want to the estimate for the second term in the right-hand side in \eqref{eq:4.1}. 
As is in Ikehata \cite{Ikehata2}, we consider indefinite integral of $v$. We define 
\begin{equation}\label{eq:w}
w(t,x) = \int_0^t v(s,x) ds.
\end{equation}
Then $w$ satisfies
\begin{equation}\label{eq:4.3}
\left\{
\begin{array}{ll}
\displaystyle ( \partial_t^2  - \triangle + B_\lambda (t,x) \partial_t ) w \\
\displaystyle = \int_0^t \left( \partial_t B_\lambda v + F_\lambda \right) ds +  B_\lambda (0) v_0 + v_1  &\quad (t,x) \in [0,\infty) \times \Omega_\lambda , \\ 
w(0,x) = 0 , \partial_t w(0,x) = v_0 (x)  &\quad x \in  \Omega_\lambda, \\
w(t,x) = 0 &\quad (t,x) \in [0,\infty) \times \partial \Omega_\lambda.
\end{array}
\right.
\end{equation}
We remark that $\partial_t w = v$ and $E(w(t))$ is well-defined in $[0,\infty)$.
\begin{lem}\label{lem:4.1}
We assume that following {\rm \textbf{(H1$)_\lambda$}} and {\rm \textbf{(H2$)_\lambda$}} hold:
\begin{description}
 \item[(H1$)_\lambda$] $\| d_0(\cdot) \{ B_\lambda(0)v_0 + v_1 \} \|_2 <\infty $,
 \item[(H2$)_\lambda$] $\displaystyle  \int_0^\infty \| d_0(\cdot) \partial_t B_\lambda (s) \|_{\infty} ds < \infty$,
\end{description}
where $d_0$ is defined in \eqref{eq:d_0}.
Then it holds following {\rm{(i)}} and {\rm(ii)}.
 \begin{description}
 \item[(i)] When $d \geq 3$, there exists a constant $E_0 = E_0(v_0,v_1)$ such that 
 \begin{equation}\label{eq:4.4}
 \int_0^t \langle v,B_\lambda v \rangle ds \leq  E_0
 \end{equation}
 \item[(ii)] When $d=2$, we assume also that {\rm \textbf{(H3$)_\lambda$}} holds.
 \item[(H3$)_\lambda$] There exists $M > 0$ such that ${\rm supp}v_0 \cup {\rm supp}v_1 \subset \{x \in \Omega_\lambda : |x| < \frac{M}{\lambda} \}$. 

Then there exists $C_{\lambda,M}>0$ such that
\begin{align}\label{eq:4.5}
 &\int_0^t \langle v,B_\lambda v \rangle ds \leq C_{\lambda,M} \|(v_0,v_1)\|^2_{H^L \times H^{L-1}} + C_{\lambda,M} \left\{ \int_0^t (1+s) Zds \right\}^2. 
\end{align}
\end{description}

\proof
Taking inner product \eqref{eq:4.3} by $\partial_t w$, we have
\begin{align*}
 & \frac{d}{dt}E(w(t)) + \langle v (t) , B_\lambda(t) v (t) \rangle  \\
 &= \langle \partial_t w (t), B_\lambda (0) v_0 + v_1 \rangle + \langle \partial_t w(t),  \int_0^t \partial_t B_\lambda (s) v(s) ds \rangle + \langle \partial_t w (t), \int_0^t F_\lambda ds \rangle.
\end{align*}
Integrating above equality over $[0,t]$, we obtain
\begin{align}\label{eq:4.6}
 &E(w(t)) + \int_0^t \langle v (s), B_\lambda v (s) \rangle ds \leq \frac{1}{2} \|v_0\|_2^2 + \langle w (t), B_\lambda(0) v_0 + v_1 \rangle \\
 & + \int_0^t \langle \partial_t w(s), \int_0^s \partial_t B_\lambda(r) v(r) dr\rangle ds + \int_0^t \langle \partial_t w(s), \int_0^s F_\lambda dr \rangle ds \nonumber\\
 &= \frac{1}{2} \|v_0\|_2^2 + (A) + (B) + (C) \nonumber.
\end{align}

First, we estimate $(A)$. Using Lemma \ref{lem:hardy}, we get
\begin{align*}
 &(A) = \langle w, B_\lambda(0) v_0 + v_1 \rangle  \leq \left\| \frac{w}{d_0(\cdot)} \right\|_2  \| d_0(\cdot) \{ B_\lambda(0) v_0 + v_1 \} \|_2 \\
 &\leq C_\lambda \| \nabla w \|_2 \| d_0(\cdot) \{ B_\lambda(0) v_0 + v_1 \} \|_2 \leq \frac{1}{4} E(w(t)) + C_\lambda \| d_0(\cdot) \{ B_\lambda(0) v_0 + v_1 \} \|_2^2.
\end{align*}
In particular, if ${\rm supp}v_0 \cup {\rm supp}v_1 \subset \{x \in \Omega_\lambda | |x| < M/\lambda \}$ then we have
\begin{equation}\nonumber
\| d(\cdot) \{ B_\lambda(0) v_0 + v_1 \} \|_2^2 \leq C_{\lambda,M} \|(v_0, v_1) \|_{H^L \times H^{L-1}}^2.
\end{equation}

Second, we estimate $(B)$. Using (H2$)_\lambda$ and \eqref{eq:prop3.2}, we calculate 
\begin{align*}
 & (B) = \int_0^t \langle \partial_t w (s), \int_0^s \partial_t B_\lambda (r) v(r) dr \rangle ds \\
 &= \langle w(t), \int_0^t \partial_t B_\lambda v(s) ds \rangle  - \int_0^t \langle w(s), \partial_t  B_\lambda(s) v(s)  \rangle ds \\
 &\leq C \sup_{0 \leq s \leq t} \left\|\frac{w(s)}{d_0(\cdot)} \right\|_2 \sup_{0 \leq s \leq t} \| v(s) \|_2 \int_0^t \| d_0(\cdot) \partial_t B_\lambda (s)  \|_\infty ds \\
 &\leq C_\lambda \sup_{0 \leq s \leq t} \| \nabla w(s) \|_2 \sup_{0 \leq s \leq t} \| v(s) \|_2 \int_0^t \| d_0(\cdot) \partial_t B_\lambda (s)  \|_\infty ds   \\
 &\leq \frac{1}{4} \sup_{0 \leq s \leq t} E(w(s)) + C_\lambda \sup_{0 \leq s \leq t} \| v(s) \|^2_2 \left\{  \int_0^t \| d_0(\cdot) \partial_t B_\lambda (s)  \|_\infty ds \right\}^2 \\
 &\leq \frac{1}{4} \sup_{0 \leq s \leq t} E(w(s)) + C_\lambda \|(v_0, v_1) \|_{H^L \times H^{L-1}}^2.
\end{align*}

Finally, we estimate $(C)$. When $d \geq 3$, using Lemma \ref{lem:GN}, we have
\begin{align*}
&\int_0^t \langle \partial_t w, \int_0^s F_\lambda dr \rangle ds = \langle w , \int_0^t F_\lambda ds \rangle - \int_0^t \langle w, F_\lambda \rangle ds  \\
&\leq 2 \sup_{0 \leq s \leq t} \| w(s) \|_{\frac{2d}{d-2}} \int_0^t \|F_\lambda \|_{\frac{2d}{d+2}} ds \leq  \frac{1}{4} \sup_{0 \leq s \leq t} E(w(s)) + C_\lambda \left\{ \int_0^t \|F_\lambda \|_{\frac{2d}{d+2}} ds \right\}^2.
\end{align*}
Using Lemma \ref{prop:sobolev}, \eqref{eq:prop3.2} and $p_l \geq 2\ (l=1,2)$, we obtain 
\begin{equation}\label{eq:4.7}
\int_0^t \|F_\lambda \|_{\frac{2d}{d+2}} ds \leq C_\lambda \int_0^t Z_0 ds \leq C_\lambda \|(v_0,v_1)\|^2_{H^{L} \times H^{L-1}}.
\end{equation}
Therefore we get
\begin{equation}\nonumber
(C) \leq \frac{1}{4}\sup_{0 \leq s \leq t}E(w(s)) + C \|(v_0,v_1) \|_{H^L \times H^{L-1}}^2.
\end{equation}
On the other hand when $d=2$, using Lemma \ref{lem:GN} and H\"{o}lder inequality, we have
\begin{align*}
 &(C) = \int_0^t \langle \partial_t w, \int_0^s F_\lambda dr \rangle ds = \langle w , \int_0^t F_\lambda ds \rangle - \int_0^t \langle w, F_\lambda \rangle ds  \\
 &\leq C_\lambda \sup_{0 \leq s \leq t} \left\|\frac{w(s)}{|\cdot|} \right\|_r \left\| \int_0^t |\cdot|F_\lambda ds \right\|_{r'} \leq C_\lambda \sup_{0 \leq s \leq t} \left\|\nabla \left( \frac{w(s)}{|\cdot|} \right) \right\|_q \int_0^t \| |\cdot| F_\lambda \|_{r'} ds,
 \end{align*}
where $r \in (2.\infty) $, $\frac{1}{r} + \frac{1}{r'} =1$ and $\frac{1}{r} = \frac{1}{q} - \frac{1}{2}$. 
From the assumption {\bf (H3$)_\lambda$} and the finite speed of propagation, it holds that 
\begin{equation}\label{eq:fsp}
{\rm supp} v(t) \cup {\rm supp} \partial_t v(t) \subset \{x \in \Omega_\lambda: |x| \leq M/\lambda+2t \} \quad (t \in [0,\infty)).
\end{equation}
Thus using \eqref{eq:fsp} and considering the same way of \eqref{eq:4.7},  we get
\begin{equation}\nonumber
\int_0^t \| |\cdot| F_\lambda \|_{r'} ds \leq C_{\lambda,M} \int_0^t (1+s) Z_0(v(s)) ds.
\end{equation}
Furthermore using Lemma \ref{lem:hardy} and H\"{o}lder inequality, we obtain
\begin{align*}
&\left\|\nabla \left( \frac{w}{|\cdot|} \right) \right\|_q \leq \left\|\frac{\nabla w}{|\cdot|} \right\|_q + \left\|\frac{w}{|\cdot|^2} \right\|_q \leq C_\lambda \|\nabla w\|_2 \left\| \frac{1}{|\cdot|} \right\|_r + C_\lambda \left\| \frac{w}{d_0} \right\|_2 \left\| \frac{d_0}{|\cdot|^2} \right\|_r \\
&\leq C_\lambda \|\nabla w \|_2 \left( \left\| \frac{1}{|\cdot|} \right\|_r + \left\| \frac{d_0}{|\cdot|^{1+\varepsilon_0}} \right\|_\infty \left\| \frac{1}{|\cdot|^{1-\varepsilon_0}} \right\|_r  \right),
\end{align*}
where $\varepsilon_0 = \frac{r-2}{2r}$.
Remember $0 \notin \Omega_\lambda$ and $r > r(1-\varepsilon_0) >2$, we get
\begin{equation}\nonumber
\|\nabla w \|_2 \left( \left\| \frac{1}{|\cdot|} \right\|_r + \left\| \frac{d_0}{|\cdot|^{1+\varepsilon_0}} \right\|_\infty \left\| \frac{1}{|\cdot|^{1-\varepsilon_0}} \right\|_r  \right) \leq C_\lambda \|\nabla w \|_2.
\end{equation}
Above estimates and Lemma \ref{lem:3.3} imply that 
\begin{align*}
 &(C) \leq C_{\lambda,M} \sup_{0 \leq s \leq t} \|\nabla w(s) \|_2 \int_0^t (1+s) Z(v(s)) ds  \\
 &\leq \frac{1}{4} \sup_{0 \leq s \leq t}E(w(s)) + C_{\lambda,M} \left\{ \int_0^t (1+s) Z(v(s)) ds \right\}^2.
\end{align*}
Combining estimates for $(A),(B),(C)$ and \eqref{eq:4.6}, we get \eqref{eq:4.4} and \eqref{eq:4.5}.
This completes the proof of Lemma \ref{lem:4.1}.
\end{lem}

\begin{rem}
If $F_\lambda$ has divergence form, we can prove {\rm (i)} if that $d=2$.
Then we do not need assume {\rm \textbf{(H3$)_\lambda$}}(See for instance: \cite{W}.).
\end{rem}

\subsection*{Proof of Theorem \ref{thm:main2}}
It is easy to see that if $\{(u_0,u_1), B\}$ satisfy $\bold{(H1)}$, $\bold{(H2)}$ and $\bold{(H3)}$, then $\{(v_0,v_1), B_\lambda\}$ satisfy $\bold{(H1)_\lambda} $, $\bold{(H2)_\lambda}$ and $\bold{(H3)_\lambda}$ respectively. 
Therefore when $d \geq 3 $, combining \eqref{eq:4.1} and \eqref{eq:4.4}, we get
\begin{equation}\label{eq:4.8}
(1+t)\{ \|v(t)\|^2_2 + Z(v(t)) \} + \int_0^t (1+s) Z(v(s)) ds \leq E_0(v_0,v_1).
\end{equation}
The above estimate means \eqref{eq:1.6}.

When $d = 2$, combining  \eqref{eq:4.1} and \eqref{eq:4.5}, we get
\begin{align*}
 &(1+t)\{ \|v(t)\|^2_2 + Z(v(t)) \} + \int_0^t (1+s) Z(v(s)) ds  \\
 &\leq C_{\lambda,M}\|(v_0,v_1)\|^2_{H^L \times H^{L-1}} + C_{\lambda,M} \left\{ \int_0^t (1+s) Z(v(s)) ds \right\}^2.
\end{align*}
The above estimate and $\displaystyle \|(v_0,v_1)\|^2_{H^L \times H^{L-1}} \leq \delta^2 $ imply
\begin{equation}\nonumber
H(t) \leq C_{\lambda,M} \delta^2 + (H(t))^2,
\end{equation}
where $H(t)= \int_0^t (1+s) Z(v(s)) ds$.
Because of $H(0) = 0$, we can choose a small $\delta$ depend on $\lambda$ and $M$ such that $H(t) \leq C_{\lambda,M}\ (t \in [0,\infty))$.
Therefore we obtain
\begin{equation}\label{eq:4.9}
(1+t)\{ \|v(t)\|^2_2 + Z(v(t)) \} + \int_0^t (1+s) Z(v(s)) ds \leq C_{\lambda,M} \|(v_0, v_1)\|^2_{H^L \times H^{L-1}} + C_{\lambda,M}^2. \nonumber
\end{equation}
This means \eqref{eq:1.6}. This completes the proof of \eqref{eq:1.6}.

Next, we prove \eqref{eq:1.7}. We calculate
\begin{equation}\nonumber
\frac{d}{dt} \left\{(1+t)^2 E(v(t)) \right\} = 2(1+t)E(v(t)) - (1+t)^2 \langle \partial_t v, B_\lambda \partial_t v \rangle + (1+t) \langle \partial_t v, F_\lambda \rangle.
\end{equation}
Integrating the above equality over $[0,t]$, we get
\begin{align*}
 &(1+t)^2 E(v(t)) + \int_0^t (1+s)^2 \langle \partial_t v, B_\lambda \partial_t v \rangle ds \\
 &\leq E(v(0)) + 2 \int_0^t (1+s) E(v(s)) ds + \int_0^t (1+s)^2 \langle \partial_t v, F_\lambda \rangle ds \\
 &\leq E(v(0)) + 2 \int_0^t (1+s) Z_0 ds + C_\lambda \left\{ \sup_{0\leq s \leq t} (1+s)^2 E(v(s)) \right\}^{\frac{1}{2}}  \int_0^t (1+s) Z_0 ds\\
 &\leq E(v(0)) + 2 \int_0^t (1+s) Z_0 ds + C_\lambda \left\{ \int_0^t (1+s) Z_0 ds \right\}^2 + \frac{1}{2} \sup_{0\leq s \leq t} (1+s)^2 E(v(s)). \nonumber
\end{align*}
Using the above estimate and \eqref{eq:1.6}, we obtain \eqref{eq:1.7}. 
This completes the proof of Theorem \ref{thm:main2}.
\qed

\appendix
\section{Estimates of nonlinear terms}
We show the estimates of the nonlinear terms $F_\lambda$.
\begin{lem}\label{lem:A1}
Let $v \in X_\delta^T$, $\delta \leq 1$ and  $|\alpha| \leq L-1$. Then there exists $C_\lambda>0$ such that 
\begin{equation}\label{eq:A1}
\|\partial^\alpha \tilde{F}_\lambda (\partial v(t)) \|_2 \leq C_{\lambda} Z(v(t)) \quad (t \in [0,T]).
\end{equation}
\proof
If $|\alpha|=0$, it is easy to see from Lemma \ref{prop:sobolev} that we prove \eqref{eq:A1}. 
Let $1 \leq |\alpha| \leq L-1$. From the chain rule, we have
\begin{align*}
 &\partial^\alpha \tilde{F}_\lambda (\partial v) \\
 &= \lambda \sum_{\alpha_1 + \cdots + \alpha_l = \alpha} C_{\alpha_1,\cdots, \alpha_l} \sum_{0 \leq j_1, \cdots, j_l \leq d \atop 1 \leq i_1, \cdots, i_l \leq d} D_{i_1j_1} \cdots D_{i_lj_l} \tilde{F}(\partial v) \partial^{\alpha_1} \partial_{j_1}v^{i_1} \cdots \partial^{\alpha_l} \partial_{j_l} v^{i_l},
\end{align*}
where $D_{ij} = D_{\xi_{ij}}$.
When $l \geq 2$, we choose $2 < q_i \leq \infty \  (i = 1,2, \cdots, l)$ satisfying
\begin{equation}\label{eq:A2}
\sum_{k=1}^l \frac{1}{q_k} = \frac{1}{2} \quad and \quad \frac{1}{2} \leq \frac{1}{q_k} + \frac{L-1-|\alpha_k|}{d}.
\end{equation}
Then it holds that 
\begin{align*}
&\|\partial^\alpha \tilde{F}_\lambda (\partial v)\|_2 \\
 &\leq C_\lambda \sum_{\alpha_1 + \cdots + \alpha_l = \alpha} \sum_{0 \leq j_1, \cdots, j_l \leq d \atop 1 \leq i_1, \cdots, i_l \leq d} \|  D_{i_1j_1} \cdots D_{i_lj_l} \tilde{F}(\partial v) \|_\infty \|\partial^{\alpha_1} \partial_{j_1}v^{i_1} \cdots \partial^{\alpha_l} \partial_{j_l} v^{i_l}\|_2   \\
 &\leq C_\lambda \sum_{\alpha_1 + \cdots + \alpha_l = \alpha} \sum_{0 \leq j_1, \cdots, j_l \leq d \atop 1 \leq i_1, \cdots, i_l \leq d} \|\partial v\|_\infty^{\max{\{0,p_1-l\}}} \prod_{k=1}^l \|\partial^{\alpha_k} \partial_{j_k} v^{i_k}\|_{q_k}\\
 &\leq C_\lambda \sum_{\alpha_1 + \cdots + \alpha_l = \alpha} \|\partial v\|_{H^{\left[\frac{d}{2} \right] + 1}}^{\max\{0,p_1-l\}} \prod_{k=1}^l \|\partial \partial^{\alpha_k} v \|_{H^{L-1-|\alpha_k|}} \\
 &\leq C_\lambda (Z(v))^{\frac{\max{\{l,p_1\}}}{2}}\leq C_\lambda Z(v),
\end{align*}
where we use the generalized H\"{o}lder inequality and a well known embedding lemma like $L^q \subset H^s\ (2<q<\infty, \frac{1}{2} \leq \frac{1}{q} + \frac{s}{d})$.
Indeed it holds from $L \geq [d/2] +3$ that 
\begin{equation}\nonumber
\sum_{k=1}^l\left\{ \frac{1}{2} - \frac{L-1-|\alpha_k|}{d} \right\} -\frac{1}{2} \leq (l-1)\left\{\frac{1}{2} + \frac{1-L}{d} \right\} \leq -\frac{l-1}{d} < 0,
\end{equation}
thus we can choose $q_k$ satisfying \eqref{eq:A2}. When $l=1$, we should choose $q_1 =2$. 
This completes the proof of Lemma \ref{lem:A1}.
\qed
\end{lem}
\begin{lem}\label{lem:A2}
Let $v \in X_\delta^T$, $\delta \leq 1$, $1 \leq |\alpha| \leq L$, $|\beta| \leq L-1$ and $|a+\beta| \leq L+1$.
For any $i,j,a$ and $b$, it holds that 
\begin{equation}\nonumber
\|\partial^\alpha v^i \partial^\beta ( c^{ab}_{ij}(\partial v) ) \|_2 \leq C_\lambda Z(v) 
\end{equation}
\proof
First, we assume $|\beta|=1$. Then we have
\begin{align*}
\|\partial^\alpha v^i \partial^\beta (c^{ab}_{ij} (\partial v)) \|_2 &\leq \| \partial^\alpha v \|_2 \sum_{0 \leq l \leq d \atop 1 \leq k \leq d} \| D_{lk} c^{ab}_{ij} (\partial v) \partial^\beta \partial_l v^k \|_\infty \\
 &\leq C \|\partial^\alpha v\|_2 \|\partial v \|_\infty^{p_2-2} \|\partial^\beta \partial v\|_\infty \leq C_\lambda (Z(v))^{\frac{p_2}{2}} \leq C_\lambda Z(v).
\end{align*}
In the same way, we can prove when $\beta=0$.

Next, we assume $|\beta| \geq 2$. In the same way as the proof of Lemma \eqref{eq:A2}, we obtain
\begin{align*}
 &\quad  \|\partial^\alpha v^i \partial^\beta (c^{ab}_{ij} (\partial v)) \|_2  \\
 &\leq C \sum_{\beta_1 + \cdots + \beta_l = \beta} \sum_{0 \leq j_1, \cdots, j_l \leq d \atop 1 \leq i_1, \cdots, i_l \leq d} \| \partial^\alpha v^i D_{i_1j_1} \cdots D_{i_lj_l} c^{ab}_{ij}(\partial v) \partial^{\beta_1} \partial_{j_1} v^{i_1} \cdots \partial^{\beta_l} \partial_{j_l}v^{i_l}  \|_2  \\
 &\leq C \sum_{\beta_1 + \cdots + \beta_l = \beta} \sum_{0 \leq j_1, \cdots, j_l \leq d \atop 1 \leq i_1, \cdots, i_l \leq d} \| D_{i_1j_1} \cdots D_{i_lj_l} c^{ab}_{ij}(\partial v) \|_\infty \| \partial^\alpha v^i \partial^{\beta_1} \partial_{j_1} v^{i_1} \cdots \partial^{\beta_l} \partial_{j_l}v^{i_l}  \|_2 \\
 &\leq C \sum_{\beta_1 + \cdots + \beta_l = \beta} \sum_{0 \leq j_1, \cdots, j_l \leq d \atop 1 \leq i_1, \cdots, i_l \leq d} \|\partial v\|_\infty^{\max{\{0,p_2-1-l\}}} \| \partial^\alpha v^i \|_{r_0} \prod_{k=1}^l \| \partial^{\beta_k} \partial_{j_k} v^{i_k}\|_{r_k} \\
 &\leq C_\lambda \sum_{\beta_1 + \cdots + \beta_l = \beta} \|\partial v\|_{H^{\left[\frac{d}{2}\right] +1}}^{\max{\{0,p_2-1-l\}}} \| \partial^\alpha v \|_{H^{L-|\alpha|}} \prod_{k=1}^l  \|\partial \partial^{\beta_k} v^{i_k}\|_{H^{L-1-|\beta_k|}} \\
 &\leq C_\lambda  Z^{\frac{\max\{l+1,p_2+1\}}{2}} \leq C_\lambda Z(v),
 \end{align*}
where we choose $2 < r_k < \infty\  (k = 0,1, \cdots, l)$ satisfying
\begin{equation}\label{eq:A3}
\sum_{k=0}^{l} \frac{1}{r_k} = \frac{1}{2},\quad \frac{1}{2} \leq  \frac{1}{r_0} + \frac{L-|\alpha|}{d} \quad {\rm and} \quad \frac{1}{2} \leq \frac{1}{r_k} + \frac{L-|\beta_k|-1}{d}.
\end{equation}
Indeed it holds from $L \geq [d/2] +3$ that 
\begin{equation}\nonumber
\frac{1}{2}-\frac{L-|\alpha|}{d} +\sum_{k=1}^l\left\{ \frac{1}{2} - \frac{L-1-|\beta_k|}{d} \right\} -\frac{1}{2} \leq \frac{l}{2} + \frac{l(1-L)}{d} + \frac{1}{d} \leq -\frac{-3l+2}{2d} < 0,
\end{equation}
thus we can choose $q_k$ satisfying \eqref{eq:A3}.
This completes the proof of Lemma \ref{lem:A2}.
\qed
\end{lem}

\section*{Acknowledgement}
I am deeply grateful to Professor Mishio Kawashita (Hiroshima University) for giving constructive comments and warm encouragement.
Furthermore I would like to express the deepest appreciation to Professor Hideo Kubo (Hokkaido University) for giving insightful suggestions.

\end{document}